\documentclass[12pt]{article}
  \usepackage{color}
\usepackage{epsfig}
\usepackage{amsfonts}
\usepackage{verbatim}
\usepackage{latexsym}
\usepackage{graphicx,amsmath}
\newcommand{\qed}{$\quad \Box$}

\makeatletter \@addtoreset{equation}{section}

\begin{document}

\newcommand{\E}{\mathbb{E}}
\newcommand{\PP}{\mathbb{P}}
\newcommand{\RR}{\mathbb{R}}
\newcommand{\I}{\mathbb{I}}

\newtheorem{theorem}{Theorem}[section]
\newtheorem{lemma}{Lemma}[section]
\newtheorem{coro}{Corollary}[section]
\newtheorem{defn}{Definition}[section]
\newtheorem{assp}{Assumption}
\newtheorem{expl}{Example}[section]
\newtheorem{remark}{Remark}[section]

\newcommand\tq{{\scriptstyle{3\over 4 }\scriptstyle}}
\newcommand\qua{{\scriptstyle{1\over 4 }\scriptstyle}}
\newcommand\hf{{\textstyle{1\over 2 }\displaystyle}}
\newcommand\hhf{{\scriptstyle{1\over 2 }\scriptstyle}}

\newcommand{\eproof}{\indent\vrule height6pt width4pt depth1pt\hfil\par\medbreak}

\def\a{\alpha} \def\d{\delta}\def\g{\gamma}
\def\e{\varepsilon} \def\z{\zeta} \def\y{\eta} \def\o{\theta}
\def\vo{\vartheta} \def\k{\kappa} \def\l{\lambda} \def\m{\mu} \def\n{\nu}
\def\x{\xi}  \def\r{\rho} \def\s{\sigma}
\def\p{\phi} \def\f{\varphi}   \def\w{\omega}
\def\q{\surd} \def\i{\bot} \def\h{\forall} \def\j{\emptyset}

\def\b{\beta} \def\de{\delta} \def\up{\upsilon} \def\eq{\equiv}
\def\ve{\vee} \def\we{\wedge}

\def\F{{\cal F}} \def\Le{{\cal L}}
\def\T{\tau} \def\G{\Gamma}  \def\D{\Delta} \def\O{\Theta} \def\L{\Lambda}
\def\X{\Xi} \def\S{\Sigma} \def\W{\Omega}
\def\M{\partial} \def\N{\nabla} \def\Ex{\exists} \def\K{\times}
\def\V{\bigvee} \def\U{\bigwedge}

\def\1{\oslash} \def\2{\oplus} \def\3{\otimes} \def\4{\ominus}
\def\5{\circ} \def\6{\odot} \def\7{\backslash} \def\8{\infty}
\def\9{\bigcap} \def\0{\bigcup} \def\+{\pm} \def\-{\mp}
\def\[{\langle} \def\]{\rangle}

\def\tl{\tilde}
\def\trace{\hbox{\rm trace}}
\def\diag{\hbox{\rm diag}}
\def\for{\quad\hbox{for }}
\def\refer{\hangindent=0.3in\hangafter=1}

\newcommand{\mz}{{m_0}}
\newcommand\wD{\widehat{\D}}

\def\bb{\begin}
\def\bc{\begin{center}}       \def\ec{\end{center}}
\def\ba{\begin{array}}        \def\ea{\end{array}}
\def\be{\begin{equation}}     \def\ee{\end{equation}}
\def\bea{\begin{eqnarray}}    \def\eea{\end{eqnarray}}
\def\beaa{\begin{eqnarray*}}  \def\eeaa{\end{eqnarray*}}
\def\la{\label}

\title{
Dynamical Behaviors of the  Tumor-immune System  in a Stochastic Environment
}
\author{Xiaoyue Li,\thanks{School of Mathematics and Statistics,
Northeast Normal University, Changchun, Jilin, 130024, China. Research
of this author was supported in part by National Natural Science Foundation of China (11171056), the Natural Science Foundation of Jilin Province (No. 20170101044JC), the Education Department of Jilin Province (No. JJKH20170904KJ).}
\and Guoting Song,\thanks{School of Mathematics and Statistics,
Northeast Normal University, Changchun, Jilin, 130024, China.}
\and Yang Xia,\thanks{School of Mathematics and Statistics,
Northeast Normal University, Changchun, Jilin, 130024, China.}
\and Chenggui Yuan\thanks{Department of Mathematics, Swansea University, Bay Campus, SA1 8EN, UK.}}
\date{}
\maketitle

\begin{abstract}

This paper investigates dynamic behaviors of the tumor-immune system perturbed by environmental noise. The model describes the response of the cytotoxic T lymphocyte~(CTL) to the growth of an immunogenic tumour. The main methods are stochastic Lyapunov  analysis, comparison theorem for stochastic differential equations (SDEs) and strong ergodicity theorem.
Firstly, we prove the existence and uniqueness of the global positive solution for the tumor-immune system. Then we go a further step to study  the boundaries of moments for tumor cells and effector cells and  the asymptotic behavior in the boundary equilibrium points. Furthermore, we discuss the existence and uniqueness of stationary distribution and stochastic permanence of the tumor-immune system. Finally,  we give several  examples and numerical simulations to verify our results.

 \vspace{3mm} \noindent {\bf Keywords.} Tumor-immune system; Stochastic permanence; Comparison theorem; Invariant measure; Ergodicity.


\end{abstract}

\newpage

\section{Introduction}\label{sec:intr}

 At present, cancer is considered to be one of the most complicated diseases to be treated clinically and one of the most dreadful killers in the world today.
 Keeping in mind its devastating nature, a great deal of human and economic resources are devoted to the research on cancer biology and subsequent development of proper therapeutic measures.  Surgery, radiation therapy, and chemotherapy are the three traditional therapy procedures that are practised for treatment of cancer.
 However, all these procedures are characterized by a relatively low efficacy and high toxicity for the patient.
Therefore, compared with traditional treatment methods, emerging immunotherapy has great development prospects. Immunotherapy, also known as biological therapy, usually refers to the use of cytokines, a protein hormone that mediate both natural and specific immunity to induce antitumor responses of immune system.

Mathematical models of tumour-immune system and their dynamical behaviors \cite{AB1997, AM2004},  help us to understand
better how host immune cells and cancerous cells evolve and interact.
 In order to get closer to reality more and more tumour-immune models have been studied, for instance, \cite{Bellomo3, DRW2005, D2005,  Natalia, MOE2016, OS1998, Villasana, Wu, Yafia} and reference therein.
    It's worth noticing that a classical mathematical simplified  tumour-immune model
   \be\la{eq1.4}
 \left\{
\begin{array}{lll}
\mathrm{d}x(t)&=&\displaystyle \left(\sigma+\frac{\rho{x(t)y(t)}}{\eta+y(t)}-\mu{x(t)y(t)}-\delta{x(t)}\right)\mathrm{d}t,\!\!\!\!\\
\mathrm{d}y(t)&=&\displaystyle \left(\alpha{y(t)}-\beta{y^{2}(t)}-x(t)y(t)\right)\mathrm{d}t \!\!\!\!
\end{array}
\right.
\ee
   is proposed to simulate the interaction of the CTL with immunogenic tumor cells and took into account the inactivation of effector cells as well as the penetration of effector cells  into tumor cells by Kuznetsov and Taylor  \cite{KMT1994},  where~$x$ represents non-dimensional local concentration  of effector cells (EC), ~$y$ represents the
non-dimensional local concentration of ~tumor cells (TC).
Their model can be applied to describe two different mechanisms of the tumor: tumor dormancy and sneaking through.
   {Yafia \cite{Yafia} studied the stability of the equilibriums and proved the existence of a family of periodic solutions bifurcating from the nontrivial steady state of the Kuznetsov-Taylor model with a delay.}
    More complete bibliography about the evolution of cells and the relevant role of cellular phenomena in directing the body toward recovery or toward illness can be found in \cite{BB2000, NM2000, PG1997}.

In the tumor tissue, the growth rate and cytotoxic parameters are influenced by many environmental factors, e.g. the degree of vascularization of tissues, the supply of oxygen, the supply of nutrients, the immunological state of the host, chemical agents, temperature, radiations, gene expression, protein synthesis and antigen shedding from the cell surface, etc. Due to the complexity, it is unavoidable that in the course of time the parameters of the system undergo random variations which give them a stochastic character \cite{EB2005,GL1978, GR1974,MAS2004}. Inclusion of randomness in mathematical models of biological and biochemical processes is thus necessary for better understanding of mechanisms which govern the biological systems. Considering the impact of the stochastic volatility of environment, we assume that environmental fluctuations mainly affect the culling rate  of effector cells ~$\delta$ and the intrinsic growth rate  of tumor cells ~$\alpha$.
$$
-\delta\mathrm{d} t\rightarrow -\delta\mathrm{d} t+\sigma_{1}\mathrm{d}B_{1}(t),~~\alpha\mathrm{d} t\rightarrow \alpha\mathrm{d} t+\sigma_{2}\mathrm{d}B_{2}(t),
$$
where $B_{1}(t)$ and $B_{2}(t)$ are the $1$-dimensional Brown motion and independent, and $\sigma_{1}$ and $\sigma_{2}$ denote the intensity of white noises. Thus  the stochastic tumor-immune model is described  by the following SDE
\be\la{eq1.5}
 \left\{\!\!\!
\begin{array}{lll}
\mathrm{d}x(t)&=&\displaystyle \left(\sigma+\frac{\rho{x(t)y(t)}}{\eta+y(t)}-\mu{x(t)y(t)}-\delta{x(t)}\right)dt+\sigma_{1}{x(t)}\mathrm{d}B_{1}(t),\!\!\!\!\\
\mathrm{d}y(t)&=&\displaystyle \left(\alpha{y(t)}-\beta{y^{2}(t)}-x(t)y(t)\right)dt+\sigma_{2}{y(t)}\mathrm{d}B_{2}(t),\!\!\!\!
\end{array}
\right.
\ee
with an initial value~$x(0)=x_{0}>0$,  $y(0)=y_{0}>0$. Based on the actual background of the model, we assume that $\sigma_{1}, \sigma_{2}$ and all other parameters are non-negative real numbers. Obviously, the model ~$(\ref{eq1.5})$ degenerates into ~$(\ref{eq1.4})$ if ~$\sigma_{1}=0$, $\sigma_{2}=0$.

In the last years, stochastic growth models for cancer cells have been  developed, one can see \cite{AG2006, FBP2000, TC2006} and reference therein. Lyapunov exponent method and Fokker-Planck method are used to investigate the stability of the stochastic models by numerical simulations.  Mukhopadhyay and Bhattacharyya \cite{MB2009} analyzed the stochastic stability  for a stochastic virus-tumor-immune model. Riccardo,  Dumitru and Oana \cite{RDO2013} studied the stochastic stability of the stochastic ~Kuznetsov-Taylor model by constructing the Lyapunov function.
  Li and Cheng \cite{LC2017} established the tumor-immune model describing the interaction and competition between the tumor cells and immune system based on the Michaelis-Menten enzyme kinetics, and gave the threshold conditions for extinction, weak persistence and stochastic persistence of tumor cells by the rigorous theoretical proofs, to name a few.

In this paper our main aim is to investigate the stochastic Kuznetsov-Taylor tumor-immune model (\ref{eq1.5}), which describes the response of the CTL to the growth of immunogenic tumor cells. Combing the stochastic Lyapunov analysis with the comparison principle for SDEs and making use of the strong ergodicity theorem, we discuss the asymptotic behaviors including the stochastic ultimately boundedness in moment, the limit distribution as well as  the ergodicity. Especially, it is pointed out that when tumor cells subject to strong stochastic perturbations, the density of tumor cells is exponentially decreasing while the density of effector cells tends to the stationary distribution. On the other hand, due to weak noises, existence and uniqueness of the  stationary distribution with the support set in $\RR^2_+$ is yielded, which implies that  tumor cells and immune cells are stochastically permanence. These obtained judgement criteria on extinction and permanence will provide us some inspirations on how to make more effective and precise therapeutic schedule to eliminate tumor cells and improve the treatment of cancer.

The rest of the paper is arranged as follow. Section 2 gives some notation and proves the existence of the unique global positive solution. Section~3 obtains the ultimate moment boundedness of the global positive solution. Section 4 yields the ergodicity of tumor cells and effector cells in the stochastic tumor-immune model which implies the stochastic permanence of cells. Section 5 presents a couple of examples and numerical simulations to illustrate our results. Section 6 concludes this paper.

\section{Global positive solution}
Throughout this paper, let $(\Omega, \F, \{\F_t\}_{t\geq 0},\PP)$ be a complete filtered probability space with $\{\F_t\}_{t\geq 0}$ satisfying the usual conditions (that is, it is right continuous and $\F_0$ contains all $\PP$-null sets). Let $B(t)=(B_1(t),
B_2(t))^T$ be an $2$-dimensional Brownian motion defined on the probability space. Let $|\cdot|$ denote both the Euclidean norm in $\RR^2$.
Also let $\RR_+=\{x\in \RR| x>0\}$ and $\RR_+^2=\{(x,y)\in \RR^2| x>0, y>0\}$. Also let  $C$ denote  a generic positive constant whose value may change in different appearances.
Moreover, let
  $C^{2,1 } (\RR^2  \times \RR_+;  {{\RR}}_+)$ denote the family of all nonnegative functions
 $V(x,  t)$ on $\RR^2 \times \RR_+$ which are continuously twice
 differentiable in $x$ and once differentiable in $t$.
 For each $V\in C^{2,1 } (\RR^2\times \RR_+;   \RR_+)$, define an
operator ${\cal{L}} $ such that ${\cal{L}}V$:
 $  \RR^2 \times \RR_+ \rightarrow \RR$ satisfying
\begin{equation}\label{eqL}\begin{split}
{\cal{L}}V(x,y, t ) &=  V_t(x,y,t)+ V_x(x,y,
t)\Big(\sigma+\frac{\rho{x y }}{\eta+y }-\mu{x y }-\delta x \Big)\\
 &~~ + V_y(x,y,
t)\Big(\alpha{y }-\beta{y^{2} }-x y \Big)+\frac{1}{2}  \Big( V_{xx}(x,y, t) \sigma^2_1x^2 +V_{yy}(x,y, t) \sigma^2_2 y^2
 \Big).
\end{split}
 \end{equation}

Since $x(t)$ represents the
 density of ~EC, ~$y(t)$ represents the
 density of ~TC, both $x(t)$ and $y(t)$ in \eqref{eq1.5} should be positive.  The theorem below gives an affirmative answer.
\begin{theorem}\label{th3.1}
For any initial value $(x_{0},y_{0})\in{\RR^{2}_{+}},$
  the equation $(\ref{eq1.5})$ has a unique global positive solution $(x(t),y(t))$ for all $t\geq 0$ with probability one.
\end{theorem}
{\bf Proof.}
Note that the coefficients of $(\ref{eq1.5})$
are locally Lipschitz continuous, so for any given initial value $(x_{0},y_{0})\in{\RR^{2}_+}$, there is a unique positive local solution $(x(t),y(t))~(t\in[0,\tau_e))$, where  $\tau_{e}$ is the explosion time . To show this solution is global, we need to prove $\tau_{e}=\infty~a.s.$
Choose  an $m_{0}\geq 1$ such that~$x_{0}\in ( {1}/{m_{0}},m_{0})$,  $y_{0}\in ( {1}/{m_{0}},m_{0})$. For any positive $m\geq  m_{0}$, define the stopping time as follows
\begin{equation}\label{eq3.1}
\tau_{m}=\inf\left\{t\in[0,\tau_{e}):\min\{x(t),y(t)\}\leq\frac{1}{m} ~\mbox{ or }~\max\{x(t),y(t)\}\geq m\right\},
\end{equation}
We set~$\inf\emptyset=\infty$, clearly, ~$\tau_{m}\leq \tau_{e}$, and~$\tau_{m}$ is increasing as~$m\rightarrow\infty$. Let~$\tau_{\infty}= \lim \limits_{m \rightarrow \infty}\tau_{m}$,  ~$\tau_{\infty}\leq\tau_{e}$ a.s. If we can prove that~$\tau_\infty=\infty~a.s.$ then ~$\tau_e=\infty~a.s.$

Here we use the proof method of contradiction.  Suppose that ~$\tau_\infty=\infty~a.s.$ doesn't  hold, then there exist constants~$T>0$ and  $\varepsilon\in(0,1)$ such that
$$\mathbb{P}(\tau_{\infty}\leq T)>\varepsilon.$$
This implies that exists an~$ m_{1}\geq m_{0}$ such that for all $m \geq m_{1}$
\begin{equation}\label{eq3.2}
\mathbb{P}(\tau_{m}\leq T)\geq \frac{\varepsilon}{2}.
\end{equation}
Define
$$
V(x,y)=(x+1-\log x)+(y+1-\log y), ~~\forall (x,y)\in {R^{2}_+}.
$$
Using the It\^{o} formula, we have
\begin{align}\label{cy1}
\mathbb{E}[V(x(\tau_{m}\wedge T),y(\tau_{m}\wedge T))] = V(x_0,y_0)+\mathbb{E}\displaystyle\int_{0}^{\tau_{m}\wedge T}\mathcal LV(x(t),y(t))\mathrm{d}t,
\end{align}
where
\begin{equation*}
\begin{split}
\mathcal LV(x,y)
=&\displaystyle (\sigma+\delta-\alpha+\frac{1}{2}\sigma_{1}^{2}+\frac{1}{2}\sigma_{2}^{2})+\frac{\rho xy}{\eta+y}+x+(\mu+\alpha+\beta)y\\
&\displaystyle-\mu xy-\delta x-\frac{\sigma}{x}-\frac{\rho y}{\eta+y}-\beta y^{2}-xy\\
\leq &\displaystyle (\sigma+\delta+\frac{1}{2}\sigma_{1}^{2}+\frac{1}{2}\sigma_{2}^{2})+(\rho+1)x+(\mu+\alpha+\beta)y.\\
\leq &\displaystyle v_{1}+2(\rho+1)(x+1-\log x)+2(\mu+\alpha+\beta)(y+1-\log y)\\
\leq &\displaystyle v_{1}+v_{2}V(x,y),
\end{split}
\end{equation*}
with ~$v_{1}=\sigma+\delta+\frac{1}{2}\sigma_{1}^{2}+\frac{1}{2}\sigma_{2}^{2}$,~~$v_{2}=2(\rho+1+\mu+\alpha+\beta)$. This, together with \eqref{cy1},
implies
$$
\mathbb{E}[V(x(\tau_{m}\wedge T),y(\tau_{m}\wedge T))] \leq V(x_0,y_0)+v_{1}{T}+v_{2}\mathbb{E}\displaystyle\int_{0}^{\tau_{m}\wedge T}V(x(t),y(t))\mathrm{d}t.
$$
The Gronwall inequality yields that
\begin{equation}\label{eq3.3}
\mathbb{E}[V(x(\tau_{m}\wedge T),y(\tau_{m}\wedge T))] \leq (V(x_0,y_0)+v_{1}T)e^{v_{2}T} .
\end{equation}
Let~$\Omega_{m}=\{\omega:\tau_{m}\leq T\}$, then~$\forall   \omega\in\Omega_{m}$, at least one of $x(\tau_{m}(\omega)\wedge T)$ and $y(\tau_{m}(\omega)\wedge T)$ is equal to~$1/m$ or ~$ {m}$.
 Hence, we have
$$
(m+1-\log m)\wedge(\frac{1}{m}+1+\log m)\leq V(x(\tau_{m}\wedge T),y(\tau_{m}\wedge T)).
$$
Due to \eqref{eq3.2} and \eqref{eq3.3}, we arrive at
\begin{equation}\label{cy2}
\begin{split}
\displaystyle\frac{\varepsilon}{2}(m+1-\log m)\wedge{(\frac{1}{m}+1+\log m)}
&\leq{\mathbb{E}[\mathbf{I}_{\Omega_{m}}(\omega)V(x(\tau_{m}\wedge {T}),y(\tau_{m}\wedge {T}))]}\\
&\leq  (V(x_0,y_0)+v_{1}T)e^{v_{2}T},
\end{split}
\end{equation}
where~$\mathbf{I}_{\Omega_{m}}(\omega)$ is  the indicate function of ~$\Omega_{m}$. On the other hand, one observes
$$
\lim\limits_{m\rightarrow\infty}(m+1-\log m)\wedge(\frac{1}{m}+1+\log m)=\infty,
$$
Taking~$m\rightarrow \infty$ in \eqref{cy2}, we obtain
$$
\infty\leq   (V(x_0,y_0)+v_{1}T)e^{v_{2}T} < \infty,
$$
which results in a contradiction. The proof is therefore complete.
\qed

\section{Moment boundedness}
Based on the existence result of positive solutions, this section focuses on the moment estimation of the processes~$x(t)$ and~$y(t)$. In order to discuss the uniform boundary of $\E [y^k(t)]$, we introduce an auxiliary process $\psi(t)$ described by
\begin{equation}\label{eq3.4}
\left\{
\begin{array}{lll}
\mathrm{d}\psi(t)=\psi(t)\left[\alpha-\beta \psi(t)\right]\mathrm{d}t+\sigma_{2}\psi(t)\mathrm{d}B_{2}(t),\\
\psi(0)=y_{0}>0,
\end{array}
\right.
\end{equation}
where~$B_{2}(t)$ is the Brownian motion defined  in~$(\ref{eq1.5})$. By utilizing a comparison theorem, one observes that ~$0<y(t)\leq \psi(t)$ for all~$t\geq0$ a.s. The following result is taken from \cite{H2}, we cite it as a lemma.
\begin{lemma}\label{le3.1}{\rm $^{\hbox{\cite{H2}}}$}
Let $\psi(t)$ be the solution of \eqref{eq3.4}. The it holds that for any~$k>1,$
\begin{equation}\label{eq3.5}
\displaystyle \mathbb{E}\psi^{k}(t)\leq  \left[\frac{1}{x_{0}}\displaystyle e^{-(\alpha+\frac{k-1}{2}\sigma^{2}_{2})t}
+\frac{2\beta}{2\alpha+(k-1)\sigma^{2}_{2}}\left(1-e^{-(\alpha+\frac{k-1}{2}\sigma^{2}_{2})t}\right)\right]^{-k}.
\end{equation}
Therefore, we have
$$
\limsup\limits_{t\rightarrow \infty } \mathbb{E}\psi^{k}(t)\leq \varrho_k:=\Big(\frac{2\alpha+(k-1)\sigma^{2}_{2}}{2\beta}\Big)^{k}, \, \forall k>1.
$$
\end{lemma}
We now investigate the asymptotic properties  of the moments of   $y(t)$.
\begin{theorem}\label{le3.2}
For any~$k>1$, we have
$$
\limsup\limits_{t\rightarrow +\infty}\mathbb{E}[y^{k}(t)]\leq \varrho_k;
$$
For any~$0<k\leq 1$, we have
$$\limsup\limits_{t\rightarrow +\infty}\mathbb{E}[y^{k}(t)]\leq
(\varrho_2)^{\frac{k}{2}}.
$$
\end{theorem}
{\bf Proof.}
Since $0<y(t)\le \psi(t)$,  If~$k>1$, by lemma~$\ref{le3.1}$, we have
$$
\limsup\limits_{t\rightarrow +\infty}\mathbb{E}[y^{k}(t)]\leq\limsup\limits_{t\rightarrow +\infty}\mathbb{E}[\psi^{k}(t)]\leq \varrho_k.
$$
Additionally,  If~$0<k\leq 1$, by~H\"{o}lder inequality, we have$$\limsup\limits_{t\rightarrow +\infty}\mathbb{E}[y^{k}(t)]\leq\limsup\limits_{t\rightarrow +\infty}\left[\mathbb{E}[y^{2}(t)]\right]^{\frac{k}{2}}\leq
(\varrho_2)^{\frac{k}{2}},$$
as required.
 \qed

Next, we continue to consider the asymptotic property of the moments of  $x(t)$. By virtue of the interaction between $x(t)$ and $y(t)$ and the positivity of $y(t)$ we provide the following sufficient result  for the moment boundedness of $x(t)$.
\begin{theorem}\label{le3.3}
 For any $\theta \in (0, 1+2\delta/\sigma_1^2)$ and  $ c>[\rho/\eta-\mu]^+$,
$$
  \limsup\limits_{t\rightarrow \infty } \mathbb{E}[(1+x(t)+cy(t))^{\theta}]\leq L(c, \theta),
$$
where $L(c, \theta)$  is a positive constant dependent on $\theta$ and $c $, which is defined by \eqref{eq+2} below.
\end{theorem}
{\bf Proof.}
Define the function~$f_1(y)=\displaystyle\frac{\rho y}{\eta+y}-(\mu+c) y$ for any $y\geq 0 $, then
$$
f_1'(y)=\displaystyle\frac{-(\mu+c)y^2-2(\mu+c)\eta y+(\rho-\eta(\mu+c))\eta}{(\eta+y)^2}.
$$
Solving $f_1'(y)=0$, we obtain two roots
$$y_1=\displaystyle-\eta-\sqrt{\frac{\eta\rho}{\mu+c}}<0,
~~~~~y_2=\displaystyle-\eta+\sqrt{\frac{\eta\rho}{\mu+c}}.$$
Therefore,
if~$y_2\leq0$, namely~$\rho\leq(\mu+c)\eta$, we then have ~
$f_1'(y)< 0$, $\forall y>0$, and ~$f_1(y)< f_1(0)=0$, $\forall y>0$.
For  a fixed ~$\theta\in (0, 1+2\delta/\sigma_1^2)$, define
$$
V_{1}(x,y)=(1+x+cy)^{\theta}, ~~\forall (x,y)\in {\RR^{2}_+}.
$$
By virtue of \eqref{eqL} computing   $\mathcal L V_1(x,y)$  leads to
\begin{align}\label{eq3.6}
\begin{split}
&\displaystyle\mathcal LV_{1}(x,y)\\
 =& \displaystyle \theta(1+x+cy)^{\theta-2}\Bigg\{\bigg[\frac{\rho y}{\eta+y}-(\mu+c)y-\delta+\frac{\theta-1}{2}\sigma^{2}_{1}\bigg]x^{2} \\
&~~~~~~~~~~~~~~~+\Big[-c\big(\beta+ \mu+c\big)y^{2}+\big(c\alpha+\frac{c\rho y}{\eta+y}-c\delta-\mu-c\big)y+\frac{\rho y}{\eta+y}-\delta+\sigma\Big]x  \\
&~~~~~~~~~~~~~~~~-c^2\beta y^{3}+\Big(c^2\alpha-c\beta+\frac{\theta-1}{2}c^2\sigma^{2}_{2}\Big)y^{2}+c(\alpha+\sigma)y+\sigma\Bigg\}.
\end{split}
\end{align}
Since $0<\theta<1+2\delta/\sigma_{1}^{2} $, that is
 ~$\displaystyle\delta+\frac{(1-\theta)}{2}\sigma_{1}^{2}>0$. Now, choose a positive constant~$\kappa:=\kappa(\theta)$ sufficiently small such that $$\displaystyle L_1(\theta):=\delta+\frac{(1-\theta)}{2}\sigma_{1}^{2}
-\frac{\kappa}{\theta}>0.$$
By the~It\^{o} formula,
\begin{equation*}
M_{v_{1}}(t):= e^{\kappa t}V_{1}(x(t),y(t))-V_{1}(x_{0},y_{0})-\displaystyle\int_{0}^{t}\mathcal L[e^{\kappa t}V_{1}(x(s),y(s))]\mathrm{d}s
\end{equation*}
is a local martingale.  And
\begin{align}\label{eq3.8}
\begin{split}
&\displaystyle\mathcal L[e^{\kappa t}V_{1}(x,y)]\\
=&\displaystyle \kappa e^{\kappa t}V_{1}(x,y)+e^{\kappa t}\mathcal LV_{1}(x,y)\\
\leq &\displaystyle\theta e^{\kappa t}(1+x+cy)^{\theta-2}
\bigg\{\Big[\frac{\rho y}{\eta+y}-(\mu+c)y-\delta+\frac{\theta-1}{2}\sigma^{2}_{1}+\frac{\kappa}{\theta}\Big]x^{2}\\
&~~~~+\Big[ -c\big(\beta+\mu +c\big)y^{2}
 +\Big(c\alpha +c\rho-c\delta-\mu-c+2c\frac{\kappa}{\theta}\Big)y
 +\rho-\delta +\sigma +\frac{2\kappa}{\theta}\Big]x\\
&~~~~ -c^2\beta y^{3}+\Big(c^2\alpha-c\beta+\frac{\theta-1}{2}c^2\sigma^{2}_{2}+c^2\frac{\kappa}{\theta}\Big)y^{2}
+c\Big(\alpha+\sigma+\frac{2\kappa}{\theta}\Big)y
+\sigma+\frac{\kappa}{\theta}\bigg\}\\
\leq  &\displaystyle \theta e^{\kappa t}(1+x+cy)^{\theta-2}(-L_{1}(\theta)x^{2}+L_{2}(c,\theta)x+L_{3}(c,\theta))\\
\leq & \displaystyle L_{4}(c,\theta)e^{\kappa t},
\end{split}
\end{align}
where
$$L_{2}(c,\theta)=\sup\limits_{y\in \mathbb{R}_{+} } \Big\{-c\big(\beta+\mu +c\big)y^{2}+\Big(c\alpha +c\rho-c\delta-\mu -c+2c\frac{\kappa}{\theta}\Big)y
+\rho-\delta +\alpha +\frac{2\kappa}{\theta}\Big\},$$
$$L_{3}(c,\theta)=\sup\limits_{y\in \mathbb{R}_{+} } \Big\{-c^2\beta y^{3}+\Big(c^2\alpha-c\beta+\frac{\theta-1}{2}c^2\sigma^{2}_{2}+c^2\frac{\kappa}{\theta}\Big)y^{2}
+c\Big(\alpha+\sigma+\frac{2\kappa}{\theta}\Big)y +\sigma+\frac{\kappa}{\theta}\Big\},$$
$$L_{4}(c,\theta)=1\vee \sup\limits_{x\in \mathbb{R}_{+} }  \big\{-L_1(\theta)x^{2}+L_{2}(c,\theta)x+L_{3}(c,\theta)\big\} .$$
 Let~$n_{0}>0$ be sufficiently large for~$x_{0},~y_0$ lying within the interval~$( {1}/{n_{0}},n_{0})$. For any constant $n\geq n_0$, define the stopping time~$$\xi_{n}=\inf\{t\geq0| ~\max\{x(t),y(t)\}\geq  n\}.$$ Note~$\xi_{n}$ is monotonically increasing and hence has a (finite or infinite) limit. Denote the limit by~$\xi_{\infty}$. For any~$n$ sufficiently large, we have~$\xi_{n}\geq \tau_{n}$, where~$\tau_{n}$ is defined by~$(\ref{eq3.1})$. By Theorem~$\ref{th3.1}$, we have ~$\tau_{\infty}=\infty$, then ~$\xi_{\infty}=\infty$. The local martingale property implies that~$\mathbb{E}[M_{v_{1}}(t\wedge\xi_{n})]=0$. That is, for any~$t\geq 0$
\begin{equation}\label{eq+1}
\mathbb{E}[e^{\kappa(t\wedge\xi_{n})}V_{1}(x(t\wedge\xi_{n}),y(t\wedge\xi_{n}))]=\mathbb{E}[V_{1}(x_{0},y_{0})]+\mathbb{E}\displaystyle\int_{0}^{t\wedge\xi_{n}}\mathcal L[e^{\kappa s}V_{1}(x(s),y(s))]\mathrm{d}s.
\end{equation}
From the definition of~$\xi_{n}$, we have ~$e^{\kappa(t\wedge\xi_{n})}(1+x(t\wedge\xi_{n})+cy(t\wedge\xi_{n}))^{\theta}$ is monotonically increasing. Let~$n\rightarrow\infty$, we obtain
$$
e^{\kappa(t\wedge\xi_{n})}(1+x(t\wedge\xi_{n})+cy(t\wedge\xi_{n}))^{\theta}~\uparrow~ e^{\kappa t}(1+x(t)+cy(t))^{\theta}~~~~a.s.
$$
By the monotone convergence theorem,
$$
\mathbb{E}[e^{\kappa(t\wedge\xi_{n})}V_{1}(x(t\wedge\xi_{n}),y(t\wedge\xi_{n}))]~\rightarrow~ \mathbb{E}[e^{\kappa t}V_{1}(x(t),y(t))],~ ~\hbox{as}~ n\rightarrow \infty.
$$
By the dominated convergence theorem,
\begin{equation*}
\displaystyle \mathbb{E}\int_{0}^{t\wedge\xi_{n}}\mathcal L[e^{\kappa s }V_{1}(x(s), y(s))]\mathrm{d}s\rightarrow \displaystyle \mathbb{E}\int_{0}^{t}\mathcal L[e^{\kappa s }V_{1}(x(s), y(s))]\mathrm{d}s,~~ \hbox{as}~ n\rightarrow \infty.
\end{equation*}
Therefore, letting~$n\rightarrow\infty$ in \eqref{eq+1} yields
\begin{equation}\label{eq3.9}
\mathbb{E}[e^{\kappa t}V_{1}(x(t),y(t))]=\mathbb{E}[V_{1}(x_{0},y_{0})]+\mathbb{E}\displaystyle\int_{0}^{t}\mathcal L[e^{\kappa s}V_{1}(x(s),y(s))]\mathrm{d}s.
\end{equation}
This implies
$$
e^{\kappa t}\mathbb{E}[(1+x(t)+cy(t))^{\theta}]\leq\mathbb{E}[(1+x_{0}+cy_{0})^{\theta}]+\frac{L_{4}(c,\theta)}{\kappa}e^{\kappa t}.
$$
Then
$$
\mathbb{E}[(1+x(t)+cy(t))^{\theta}]\leq\mathbb{E}[(1+x_{0}+cy_{0})^{\theta}]e^{-\kappa t}+\frac{L_{4}(c,\theta)}{\kappa}.
$$
Letting~$t\rightarrow\infty$, we have
\begin{equation}\label{eq+2}
\limsup\limits_{t\rightarrow \infty } \mathbb{E}[(1+x(t)+cy(t))^{\theta}]\leq\frac{L_{4}(c,\theta)}{\kappa}=: L(c,\theta).
\end{equation}
The proof is complete.
 \qed

The positivity of $y(t)$ implies the follow result directly.
\begin{coro}\label{c1}
For any $\theta \in (0, 1+2\delta/\sigma_1^2)$ and  $ c>[\rho/\eta-\mu]^+$,
$$
  \limsup\limits_{t\rightarrow \infty } \mathbb{E}[(1+x(t) )^{\theta}]\leq L(c, \theta),
$$
where $L(c, \theta)$  is defined in Theorem \ref{le3.3}.
 \end{coro}

Due to the inequality direction in stochastic analysis it is difficult to find the lower bound of the moment of $x(t)$. Alternatively, we try to  look for the upper bound of the moment of  $1/x(t) $. Thus we get the following result.
\begin{lemma}\label{le3.4}
If ~$\theta\in (0, 2)$,  then  there exists an $L>0$ such that
\begin{equation*}
\limsup\limits_{t\rightarrow \infty } \mathbb{E}[{ {x}}^{-\theta}(t)]\leq L.
\end{equation*}
\end{lemma}
{\bf Proof.}
Let
\begin{equation*}
V_{2}({x})=\left(1+\frac{1}{x}\right)^{\theta}, ~~\forall {x}>0.
\end{equation*}
Choosing a positive constant~$\kappa$ and applying the~It\^{o} formula lead to
\begin{equation}\label{eq3.15}
M_{v_{2}}(t):= e^{\kappa t}V_{2}(  {x}(t))-V_{2}( {x}_{0})-\displaystyle\int_{0}^{t}\mathcal L_{x}[e^{\kappa s}V_{2}( {x}(s))]\mathrm{d}s
\end{equation}
is a local martingale, where
\begin{align*}
\begin{split}
&\mathcal L_x [e^{\kappa t}V_{2}( {x})]\\
:=&\displaystyle \theta e^{\kappa t}\left(1+ \frac{1}{x}\right)^{\theta-2}\bigg[-\frac{\sigma} {{x}^{3}}-\Big(\sigma-\delta
-\frac{\theta+1}{2}\sigma_{1}^{2}-\frac{\kappa}{\theta}\Big) \frac{1}{{x}^{2}}
+\Big(\delta+\sigma_{1}^{2}+\frac{2\kappa}{\theta}\Big)\frac{1}{x} \\
&\displaystyle~~~~~~~~~~~~~~~~~~~~~~
-\frac{\rho  y}{x(\eta+y)}
+\frac{\mu y}{x}-\frac{\rho y}{ {x}^{2}(\eta+y)}+\mu \frac{y}{{x}^{2}}+\frac{\kappa}{\theta}\bigg]
\end{split}
\end{align*}
Using the Young inequality yields
\begin{eqnarray} \label{eq3.16}
&&\mathcal L_x [e^{\kappa t}V_{2}( {x})]\nonumber\\
&\leq &\displaystyle \theta e^{\kappa t}\left(1+ \frac{1}{x}\right)^{\theta-2}\bigg[-\frac{\sigma} {{x}^{3}}-\Big(\sigma-\delta
-\frac{\theta+1}{2}\sigma_{1}^{2}-\frac{\kappa}{\theta}\Big) \frac{1}{{x}^{2}}
+\Big(\delta+\sigma_{1}^{2}+\frac{2\kappa}{\theta}\Big)\frac{1}{x}\nonumber \\
&&\displaystyle~~~~~~~~~~~~~~~~~~~~~~
+\frac{\mu y}{x}+\mu \frac{y}{{x}^{2}}+\frac{\kappa}{\theta}\bigg]\nonumber\\
&\leq &\displaystyle \theta e^{\kappa t}\left(1+ \frac{1}{x}\right)^{\theta-2}\bigg[-\frac{\sigma} {{x}^{3}}+\frac{4\mu }{5{x}^{\frac{5}{2}}}
-\Big(\sigma-\delta
-\frac{\theta+1}{2}\sigma_{1}^{2}-\frac{\kappa}{\theta}-\frac{\mu}{2}\Big) \frac{1}{{x}^{2}}
\nonumber\\
&&\displaystyle~~~~~~~~~
+\Big(\delta+\sigma_{1}^{2}+\frac{2\kappa}{\theta}\Big)\frac{1}{x}
+\frac{\kappa}{\theta}+\frac{\mu}{5}y^{5}+\frac{\mu}{2}y^{2}\bigg]\nonumber\\
&\leq &\displaystyle e^{\kappa t}\Big(L_{6}+\frac{\mu}{5}y^{5}+\frac{\mu}{2}y^{2}\Big),
 \end{eqnarray}
where
$$L_{6}=\sup\limits_{ {x}\in \mathbb{R}_{+} } \Big\{-\sigma {x}^{-3}+\frac{4}{5}\mu {x}^{-\frac{5}{2}}-\Big(\sigma-\delta-\frac{\theta+1}{2}\sigma_{1}^{2}
-\frac{\kappa}{\theta}-\frac{\mu}{2}\Big) {x}^{-2}+\Big(\delta+\sigma_{1}^{2}+\frac{2\kappa}{\theta}\Big) {x}^{-1}+\frac{\kappa}{\theta}\Big\}.$$
Let~$n_{0}>0$ be sufficiently large for the initial value  ~$ {x}_{0}$  lying within the interval~$({1}/{n_{0}}, n_{0})$. For any $n\geq n_0$, define the stopping time~$$\tilde{\xi}_{n}=\inf\{t\geq0,  {x}(t)\leq 1/ n\}. $$ Note~$\tilde{\xi}_{n}$ is monotonically increasing and hence has a (finite or infinite) limit. Denote the limit by~$\tilde{\xi}_{\infty}$. For any~$n$ sufficiently large, we have~$\tilde{\xi}_{n}\geq \tau_{n}$, where~$\tau_{n}$ is defined by $(\ref{eq3.1})$. By Theorem~$\ref{th3.1}$, we have~$\tau_{\infty}=\infty$, so ~$\tilde{\xi}_{\infty}=\infty$. The local martingale property implies that~$\mathbb{E}[M_{v_{2}}(t\wedge\tilde{\xi}_{n})]=0$. That is, for any~$t\geq 0$
$$
\mathbb{E}[e^{\kappa(t\wedge\tilde{\xi}_{n})}V_{2}( {x}(t\wedge\tilde{\xi}_{n}))]=\mathbb{E}[V_{2}( {x}_0)]+\mathbb{E}\displaystyle\int_{0}^{t\wedge\tilde{\xi}_{n}}\mathcal L_x[e^{\kappa s}V_{2}( {x}(s))]\mathrm{d}s.
$$
From the definition of~$\tilde{\xi}_{n}$, we have ~$e^{\kappa(t\wedge\tilde{\xi})}(1+\frac{1}{x(t\wedge\tilde{\xi}_{n})})^{\theta}$ is monotonically increasing. Letting~$n\rightarrow\infty$ yields
$$
e^{\kappa(t\wedge\tilde{\xi}_{n})}\left(1+\frac{1}{x (t\wedge\tilde{\xi}_{n})}\right)^{\theta}~\uparrow~ e^{\kappa t}\left(1+\frac{1}{x (t)}\right)^{\theta}~~~~a.s.
$$
By the monotone convergence theorem one notices that as $n\rightarrow \infty$
$$
\mathbb{E}[e^{\kappa(t\wedge\tilde{\xi}_{n})}V_{2}( {x}(t\wedge\tilde{\xi}_{n}))]~\rightarrow~ \mathbb{E}[e^{\kappa t}V_{2}( {x}(t))].
$$
Noting that ~$\mathbb{E}[y^{5}(t)]$ and  $\mathbb{E}\left[y^{2}(t)\right]$ are   bounded uniformly with respect to~$t\in(0, \infty)$, by the~Fubini theorem and \eqref{eq3.16}, we obtain
\begin{equation*}
\begin{array}{lll}
\displaystyle \mathbb{E}\int_{0}^{t}\mathcal L[e^{\kappa s }V_{2}({x}(s))]\mathrm{d}s&\leq &\displaystyle \mathbb{E}\int_{0}^{t} e^{\kappa s}\Big(L_{6}+\frac{\mu}{5}y^{5}(s)+\frac{\mu}{2}y^{2}(s)\Big)\mathrm{d}s\\
&=&\mathbb{E}\displaystyle\int_{0}^{t}L_{6}e^{\kappa s}\mathrm{d}s+\displaystyle\int_{0}^{t}e^{\kappa s}\left[\frac{\mu}{5}\mathbb{E}\big(y^5(s)\big)
+\frac{\mu}{2}\mathbb{E}\big(y^2(s)\big)\right]\mathrm{d}s<\infty.
\end{array}
\end{equation*}
Using the dominated convergence theorem implies that as $n\rightarrow \infty$
\begin{equation}
\displaystyle \mathbb{E}\int_{0}^{t\wedge\tilde{\xi}_{n}}\mathcal L[e^{\kappa s }V_{2}({x}(s))]\mathrm{d}s\rightarrow \displaystyle \mathbb{E}\int_{0}^{t}\mathcal L[e^{\kappa s }V_{2}({x}(s))]\mathrm{d}s.
\end{equation}
Therefore,  letting~$n\rightarrow\infty$ yields
\begin{equation}\label{eq3.17}
\mathbb{E}[e^{\kappa t}V_{2}( {x}(t))]=\mathbb{E}[V_{2}(x_0)]+\mathbb{E}\displaystyle\int_{0}^{t}\mathcal L[e^{\kappa s}V_{2}( {x}(s))]\mathrm{d}s.
\end{equation}
This together with Theorem \ref{le3.2} implies
\begin{align}
\begin{split}
e^{\kappa t}\mathbb{E}\left[\left(1+\frac{1}{{x}(t)}\right)^{\theta}\right]
\leq &\displaystyle \left(1+\frac{1}{{x}_0}\right)^{\theta}+\mathbb{E}\displaystyle\int_{0}^{t}e^{\kappa s}\Big(L_{6}+\frac{\mu}{5}y^{5}(s)+\frac{\mu}{2}y^{2}(s)\Big)\mathrm{d}s\\
\leq &\displaystyle \left(1+\frac{1}{{x}_0}\right)^{\theta}+\displaystyle\int_{0}^{t}e^{\kappa s}\Big[L_{6}+\frac{\mu}{5}\mathbb{E}\big(y^{5}(s)\big)+\frac{\mu}{2}\mathbb{E}\big(y^{2}(s)\big)\Big]\mathrm{d}s\\
\leq &\displaystyle \left(1+\frac{1}{{x}_0}\right)^{\theta}+\frac{L_{7}}{\kappa}e^{\kappa t},
\end{split}
\end{align}
where~$\displaystyle L_{7}:=L_{6}+\frac{\mu}{5}\sup_{t\geq 0}\mathbb{E}\big(y^{5}(t)\big)+\frac{\mu}{2}\sup_{t\geq 0}\mathbb{E}\big(y^{2}(t)\big)<\infty$, hence
$$
\mathbb{E}\left[\left(1+\frac{1}{{x}(t)}\right)^{\theta}\right]\leq\left(1+\frac{1}{{x}_0}\right)^{\theta}e^{-\kappa t}+\frac{L_{7}}{\kappa}.
$$
 We therefore obtain
$$
\limsup\limits_{t\rightarrow \infty } \mathbb{E}[ {x}^{-\theta}(t)]\leq \limsup\limits_{t\rightarrow \infty }\mathbb{E}\left[\left(1+\frac{1}{{x}(t)}\right)^{\theta}\right]\leq\frac{L_{7}}{\kappa}=: L.
$$
The proof is complete.
\qed

\section{Existence and  uniqueness of invariant measure }
This section is devoted to analyze the invariant measure.
Define~function $$f_2(y)=\displaystyle\frac{\rho y}{\eta+y}-\mu y,~~~y\geq0.$$
 Similar to the analysis of the function~$f_1(y)$ in Theorem \ref{le3.3}, we obtain
\begin{itemize}
\item[$(\mathrm{i})$]   If~ $\rho\leq\mu\eta$, ~$f_2(y)<  0$, ~$\forall y> 0$.
\item[$(\mathrm{ii})$]    If~$\rho>\mu\eta$, ~$f_2(y)\leq  (\sqrt{\rho}-\sqrt{\mu\eta})^2,$ ~$\forall y> 0$.
\end{itemize}
This implies that  for any~$y>0$,  $f_2(y)\leq[(\sqrt{\rho}-\sqrt{\mu\eta})\vee0]^2.$
We now introduce a new auxiliary process $\varphi(t)$ with respect to  $x(t)$  described by
\begin{equation}\label{eq3.18}
\begin{cases}{}
\mathrm{d}\varphi(t)&=[\sigma-(\delta-h^2)\varphi(t)]\mathrm{d}t+\sigma_{1}\varphi(t)\mathrm{d}B_{1}(t),\\
\varphi(0)&=x_{0}>0,
\end{cases}
\end{equation}
where
~\begin{equation}\label{eq-f}
 h:=(\sqrt{\rho}-\sqrt{\mu\eta})\vee0.
 \end{equation}
  If~$\delta-h^2>0$, by solving the Fokker-Planck equation (see details in \cite{DND2015}), the process~$\varphi(t)$ has a unique stationary distribution~$\nu(\cdot)$  which is the inverse Gamma distribution with parameter
$$a_1=\displaystyle\frac{2(\delta-h^2)}{\sigma_{1}^{2}}+1,~
~~b_1=\displaystyle\frac{2\sigma}{\sigma_{1}^{2}},$$  with a notation abuse slightly, we write $\phi \sim IG(a_1,~b_1)$, with probability density
$$
\displaystyle f^{*}(x)=\frac{b_1^{a_1}}{\Gamma(a_1)}x^{-(a_1+1)}e^{-\frac{b_1}{x}},~~~~~ x>0.
$$
For any~$p>0$, by the strong law of large numbers
we deduce that
\begin{equation}\label{eq3.19}
\displaystyle\lim\limits_{t\rightarrow \infty }\frac{1}{t}\int_{0}^{t}\varphi^{p}(s)\mathrm{d}s=\int_{0}^{\infty}x^{p}f^{*}(x)\mathrm{d}x:= M_{p}~~~~~a.s.
\end{equation}
Especially, if~$p=1$, ~$M_{1}=\displaystyle\frac{\sigma}{\delta-h^2}$. Moreover, ~the stationary distribution of $\displaystyle\frac{1}{\varphi(t)}$ is the Gamma distribution with parameter~$a_1$ and~$b_1$, see details in \cite{DHY2015}.
Therefore, by the~It\^{o} formula and the strong law of large numbers, noting that the mean of Gamma distribution is $a_1/b_1,$  we arrive at
\begin{equation}\label{eq3.20}
\displaystyle\lim\limits_{t\rightarrow \infty }\frac{1}{t}\ln\varphi(t)=\lim\limits_{t\rightarrow \infty }\frac{1}{t}\int_{0}^{t}\left(\frac{\sigma}{\varphi(s)}-\delta+h^2-\frac{\sigma_{1}^{2}}{2}\right)ds+\sigma_{1}\lim\limits_{t\rightarrow \infty }\frac{B_{1}(t)}{t}=0,
\end{equation}
By virtue of the comparison theorem it follows that ~$0<x(t)\leq \varphi(t)$ for all~$t\geq0$ a.s.
This implies,
\begin{equation}\label{eq3.38}
\displaystyle\limsup\limits_{t\rightarrow \infty }\frac{1}{t}\ln x(t)\leq0 ~~~~\hbox{a.s.}
\end{equation}
Furthermore, we derive the following result from \eqref{eq3.19}.
\begin{lemma}\label{le3.5}
If~$\delta-h^2>0$,
\begin{equation}\label{eq3.35}
\limsup\limits_{t\rightarrow\infty}\frac{1}{t}\int_{0}^{t}x^{p}(s)\mathrm{d}s\leq M_{p},~~\forall p>0,~~~~~a.s.
\end{equation}
Moreover,
$$
\limsup\limits_{t\rightarrow\infty}\frac{1}{t}\int_{0}^{t}x(s)\mathrm{d}s\leq \frac{\sigma}{\delta-h^2}~~~~~\hbox{a.s.}
$$
\end{lemma}

Now, we consider the auxiliary process ~$\psi(t)$ defined by \eqref{eq3.4}. If~$\displaystyle2\alpha< \sigma_{2}^{2}$, it can be easily verified that~$\lim\limits_{t\rightarrow\infty}\psi(t)=0~~a.s$. If~$\displaystyle 2\alpha>\sigma_{2}^{2}$,
by solving the Fokker-{Planck equation (see details in \cite{DHY2015}), the process~$\psi(t)$ has a unique stationary distribution~$\lambda(\cdot)$, and obeys the Gamma distribution with parameter
~$$
a_2=\displaystyle\frac{2\alpha}{\sigma_{2}^{2}}-1,
~~~~b_2=\displaystyle\frac{2\beta}{\sigma_{2}^{2}},
$$  with a notation abuse slightly, we write $\psi \sim G (a_2,~b_2)$, with density$$
\displaystyle g^{*}(x)=\frac{(b_2)^{a_2}}{\Gamma(a_2)}x^{a_2-1}e^{-b_2 x},~~~  x>0.
$$
For any~$p>0$, by the strong law of large numbers  we derive that
\begin{equation}\label{eq3.21}
\displaystyle\lim\limits_{t\rightarrow \infty }\frac{1}{t}\int_{0}^{t}\psi^{p}(s)\mathrm{d}s=\int_{0}^{\infty}x^{p}g^{*}(x)\mathrm{d}x:= \bar{M}_{p}~~~~~a.s.
\end{equation}
In particular, if~$p=1$, we have~$\displaystyle \bar{M}_{1}=\frac{1}{\beta}(\alpha-\frac{\sigma_{2}^{2}}{2})$.
Therefore, using the~It\^{o} formula implies
\begin{equation}\label{eq3.22}
\displaystyle\lim\limits_{t\rightarrow \infty }\frac{1}{t}\ln\psi(t)=\lim\limits_{t\rightarrow \infty }\frac{1}{t}\int_{0}^{t}\left(\alpha-\frac{\sigma_{2}^{2}}{2}-\beta\psi(s)\right)\mathrm{d}s+\sigma_{2}\lim\limits_{t\rightarrow \infty }\frac{B_{2}(t)}{t}=0.
\end{equation}
By virtue of the comparison theorem it follows that ~$0<y(t)\leq \psi(t)$ for all~$t\geq0$ a.s. One observes that
\begin{equation}\label{eq3.39}
\displaystyle\limsup\limits_{t\rightarrow \infty }\frac{1}{t}\ln y(t)\leq0~~~~~a.s.
\end{equation}
Furthermore, we yield the following result from \eqref{eq3.21}.
\begin{lemma}\label{le3.6}
If~$\displaystyle2\alpha> \sigma_{2}^{2}$,
\begin{equation}\label{eq3.36}
\limsup\limits_{t\rightarrow\infty}\frac{1}{t}\int_{0}^{t}y^{p}(s)\mathrm{d}s\leq \bar{M}_{p},~~\forall p>0,~~~~~a.s.
\end{equation}
Moreover,
$$
\limsup\limits_{t\rightarrow\infty}\frac{1}{t}\int_{0}^{t}y(s)\mathrm{d}s\leq \frac{1}{\beta}(\alpha- \frac{\sigma_{2}^{2}}{2})~~~~~a.s.
$$
\end{lemma}

To obtain more properties of the solution, we go a further step to consider the equation on the boundary
\begin{equation}\label{eq3.23}
\left\{
\begin{array}{l}
\mathrm{d}z(t)=(\sigma-\delta z(t))\mathrm{d}t+\sigma_{1}z(t)\mathrm{d}B_{1}(t),~~~~\forall t\geq t_0,\\
z(t_0)=x(t_0) >0,
\end{array}
\right.
\end{equation}
where ~$t_0\geq 0$  will be chosen latter. By solving the Fokker-Planck equation (see details in \cite{DND2015}), the process~$z(t)$ has a unique stationary distribution~$\mu(\cdot)$, and obeys the inverse Gamma distribution with parameter $$a_3=\displaystyle\frac{2\delta}{\sigma_{1}^{2}}+1,~~
b_3=\displaystyle\frac{2\sigma}{\sigma_{1}^{2}}.$$  With a notation abuse slightly, we write $  z \sim IG(a_3,~b_3),$  with probability density
$$
\displaystyle p^{*}(x)=\frac{ (b_3 )^{a_3}}{\Gamma(a_3)}x^{-(a_3+1)}e^{-\frac{b_3}{x}},~~~~~ x>0.
$$

In the following, we will reveal the long-time behavior  of the tumor cells and the effector cells,
 if the intensity of the noise $\sigma_2$ is large sufficiently.
\begin{theorem}\label{th3.5}
If~$\lambda_1:=\displaystyle\frac{\sigma_{2}^{2}}{2}-\alpha>0$, then we have  $$
\displaystyle \limsup\limits_{t\rightarrow\infty}\frac{\ln y(t)}{t}\leq -\lambda_1,
$$ and the distribution of ~$x(t)$ converges weakly to a unique invariant probability measure~$\pi_1(\cdot)$.
\end{theorem}
 {\bf Proof.}
 If~$\displaystyle\alpha<\frac{\sigma_{2}^{2}}{2}$, by virtue of the~It\^{o} formula, it follows from  \eqref{eq3.4} that
$$
\displaystyle\limsup\limits_{t\rightarrow \infty }\frac{1}{t}\ln\psi(t)\leq -\lambda_1~~~~~a.s.
$$
which implies that
\begin{equation}\label{eq3.24}
\displaystyle\limsup\limits_{t\rightarrow \infty }\frac{1}{t}\ln y(t)\leq -\lambda_1~~~~~a.s.
\end{equation}
For any~$\varepsilon>0$, let $t_{0}>0$ be sufficiently large such that ~$\mathbb{P}(\Omega_{\varepsilon})>1-\varepsilon$, where
$$
\displaystyle\Omega_{\varepsilon}:=\displaystyle\left\{y(t)\leq\exp\left(-\frac{\lambda_1t}{2}\right),~\forall t\geq t_{0}\right\}=\displaystyle\left\{\ln y(t)\leq-\frac{\lambda_1t}{2},~\forall t\geq t_{0}\right\},
$$  and ~$$\displaystyle\max\left\{ \frac{2\mu}{\lambda_1}\exp\left(-\frac{\lambda_1t_{0}}{2}\right), \frac{2\rho}{\lambda_1\eta}\exp\left(-\frac{\lambda_1t_{0}}{2}\right)\right\}<\frac{\varepsilon}{2},$$

\indent   {\bf Case~$(1)$.}  If~$\rho\leq\mu\eta$, ~$\displaystyle f_2(y)=\frac{\rho y}{\eta+y}-\mu y\leq0$, by the comparison theorem, we have~$\mathbb{P}\{z(t)\geq x(t),~~\forall t\geq t_{0}\}=1$. By the~It\^{o} formula, we deduce that  for almost all~$\omega\in\Omega_\varepsilon$, $\forall t\geq t_0$,
\begin{align*}
\displaystyle 0\leq \ln z(t)-\ln x(t)=&\displaystyle\sigma\int^{t}_{t_{0}}\left(\frac{1}{z(s)}-\frac{1}{x(s)}\right)\mathrm{d}s-\int^{t}_{t_{0}}\frac{\rho y(s)}{\eta+y(s)}\mathrm{d}s+\mu\int^{t}_{t_{0}}y(s)\mathrm{d}s\\
\leq&\displaystyle\mu\int^{t}_{t_{0}}\exp\left(-\frac{\lambda_1s}{2}\right)\mathrm{d}s\\
=&\displaystyle \frac{2\mu}{\lambda_1}\left[\exp\left(-\frac{\lambda_1t_{0}}{2}\right)-\exp\left(-\frac{\lambda_1t}{2}\right)\right]<\displaystyle\frac{\varepsilon}{2}.
\end{align*}
\indent   {\bf Case~$(2)$.}  If~$\rho>\mu\eta$, ~$f_2(0)>0$. Due to \eqref{eq3.24} and the continuity of $f_2(y)$ at $y=0$, one may choose $t_1\geq t_0$ such that for all $t\geq t_1$, $y(t)$ is  sufficiently small such that $\displaystyle f_2(y(t))=\frac{\rho y(t)}{\eta+y(t)}-\mu y(t)>0$. By the comparison theorem, we have~$\mathbb{P}\{x(t)\geq z(t),~~\forall t\geq t_{1}\}=1$. By the~It\^{o} formula we deduce that, for almost all~$\omega\in\Omega_\varepsilon$, $\forall t\geq t_1$,
\begin{align*}
\displaystyle 0\leq \ln x(t)-\ln z(t)=&\displaystyle\sigma\int^{t}_{t_{1}}
\left(\frac{1}{x(s)}-\frac{1}{z(s)}\right)\mathrm{d}s
+\int^{t}_{t_{1}}\frac{\rho y(s)}{\eta+y(s)}\mathrm{d}s-\mu\int^{t}_{t_{1}}y(s)\mathrm{d}s\\
\leq & \frac{\rho}{\eta} \int^{t}_{t_{1}}  y(s) \mathrm{d}s
\leq \displaystyle\frac{\rho}{\eta}\int^{t}_{t_{1}}
\exp\left(-\frac{\lambda_1s}{2}\right)\mathrm{d}s\\
=&\displaystyle \frac{2\rho}{\lambda_1\eta}\left[\exp\left(-\frac{\lambda_1t_{1}}{2}\right)-\exp\left(-\frac{\lambda_1t}{2}\right)\right]<\displaystyle\frac{\varepsilon}{2}.
\end{align*}
Therefore
\begin{equation}\label{eq3.25}
\mathbb{P}\Big\{\big|\ln z(t)-\ln x(t)\big|> \varepsilon\Big\}\leq 1-\mathbb{P}(\Omega _{\varepsilon})<\varepsilon,~~~~\forall t\geq t_{1}.
\end{equation}
Let~$\pi_1^*(\cdot)$
 be the invariant measure of~$\ln z(t)$.  In order to  show  that the distribution of~$x(t)$ converges weakly to a probability measure~$\pi_{1}(\cdot)$,
  we only need to prove that the distribution of~$\ln x(t)$ converges weakly to~$\pi_1^*(\cdot)$.
  Let
  $\mathcal\mathbb{P}(\mathbb{R})$
  represents the family of all probability measures  on~$\mathbb{R}$.
  For any $\mathbb{P}_{1}, \mathbb{P}_{2}\in{\mathbb{P}(\mathbb{R})}$,
   define the  distance as in \cite{AA}
$$
\mathrm{d}_{\mathbb{L}}(\mathbb{P}_{1},~\mathbb{P}_{2})=\sup\limits_
{f\in{\mathbb{L}}}\Big|\int_{\mathbb{R}}f(x)\mathbb{P}_{1}(\mathrm{d}x)-  \int_{\mathbb{R}}f(x)\mathbb{P}_{2}(\mathrm{d}x)\Big|,
$$
where
$$
\mathbb{L}=\Big\{f: \mathbb{R}\rightarrow \mathbb{R}:|f(x)-f(y)|\leq|x-y|~\text{ and }~|f(\cdot)|\leq 1\Big\}.
$$
By the Portmanteau theorem, we need to prove that for any~$f \in \mathbb{L},$
$$
\mathbb{E}f(\ln x(t))\rightarrow \bar{f}:=\int_{\mathbb{R}} f(x)\pi_1^*(\mathrm{d}x)=\int^{\infty}_{0}f(\ln x)\pi_1(\mathrm{d}x).
$$
Since~
the diffusion is nondegenerate, it is well known that as~$t\rightarrow\infty$ the distribution of~$\ln z(t)$ converges weakly to the unique stationary distribution~$\pi_1^*(\cdot)$, namely,
\begin{equation}\label{eq3.26}
\lim\limits_{t\rightarrow\infty} \mathbb{E}f(\ln z(t)) =\bar{f}.
\end{equation}
We now compute
\begin{align}\label{eq3.27}
\begin{split}
&\big|\mathbb{E}f(\ln x(t))-\bar{f}\big|\\
\leq&\big|\mathbb{E}f(\ln x(t))-\mathbb{E}f(\ln z(t))\big|+\big|\mathbb{E}f(\ln z(t))-\bar{f}\big|\\
=&\big|\mathbb{E}[f(\ln x(t))-f(\ln z(t))]\big|+\big|\mathbb{E}f(\ln z(t))-\bar{f}\big|\\
\le & \mathbb{E}[ |f(\ln x(t))-f(\ln z(t))|\mathbf{I}_{\{|\ln x(t) -  \ln z(t) |\leq\varepsilon\}}] \\
&+ \mathbb{E}[|f(\ln x(t))-f(\ln z(t))|\mathbf{I}_{\{|\ln x(t) -  \ln z(t) | >\varepsilon\}}]
 +\big|\mathbb{E}f(\ln z(t))-\bar{f}\big|\\
\leq& \varepsilon\mathbb{E}[\mathbf{I}_{\{|\ln x(t) -  \ln z(t) |\leq\varepsilon\}}] + 2\mathbb{E}[\mathbf{I}_{\{|\ln x(t) -
\ln z(t) |>\varepsilon}\}]
 +\big|\mathbb{E}f(\ln z(t))-\bar{f}\big|\\
=&\varepsilon
+2\mathbb{P}\Big\{\big|\ln z(t)-\ln x(t)\big|>\varepsilon\Big\}
 \!+\big|\mathbb{E}f(\ln z(t))-\bar{f}\big|.
\end{split}
\end{align}
This, together with ~$(\ref{eq3.25})$ and $(\ref{eq3.26})$, yields
$$
\limsup\limits_{t\rightarrow\infty}|\mathbb{E}f(\ln {x}(t))-\overline{f}|=0.
$$
The proof is therefore complete.
\qed

In order to investigate the probability law for the small noises we prepare two lemmas.
\begin{lemma}\label{le3.7}
 If~$\displaystyle\frac{\sigma_{2}^{2}}{2}<\alpha$,  then the property
\begin{equation}\label{eq3.28}
\displaystyle \limsup\limits_{t\rightarrow\infty}\frac{1}{t}\int^{t}_{0}\frac{1}{x(s)}\mathrm{d}s\leq \lambda_{2}
\end{equation}
holds, where~$\lambda_{2}:= \displaystyle \frac{1}{\sigma}\Big[\frac{\mu}{\beta}\Big(\alpha-\frac{\sigma_{2}^{2}}{2}\Big)
+\delta+\frac{\sigma_{1}^{2}}{2}\Big].$
\end{lemma}
{\bf Proof.}
For any~$(x_0,y_0)\in\mathbb{R}_{+}^{2}$, using the fact  $0<y(t)\leq \psi(t)$ a.s. and the  It\^{o} formula, we have
$$
\begin{array}{lll}
\displaystyle \frac{1}{t}\ln x(t)&=&\displaystyle \frac{1}{t}\int^{t}_{0}\left(\frac{\sigma}{x(s)}+\frac{\rho y(s)}{\eta+y(s)}-\mu y(s)\right)\mathrm{d}s-\delta-\frac{\sigma_{1}^{2}}{2}+\frac{\ln x_{0}}{t}+\frac{\sigma_{1}B_{1}(t)}{t}\\
&\geq&\displaystyle \frac{1}{t}\int^{t}_{0}\frac{\sigma}{x(s)}\mathrm{d}s-\frac{1}{t}\int^{t}_{0}\mu y(s)\mathrm{d}s-\delta-\frac{\sigma_{1}^{2}}{2}+\frac{\ln x_{0}}{t}+\frac{\sigma_{1}B_{1}(t)}{t}\\
&\geq&\displaystyle \frac{1}{t}\int^{t}_{0}\frac{\sigma}{x(s)}\mathrm{d}s-\frac{1}{t}\int^{t}_{0}\mu\psi(s)\mathrm{d}s-\delta-\frac{\sigma_{1}^{2}}{2}+\frac{\ln x_{0}}{t}+\frac{\sigma_{1}B_{1}(t)}{t}.
\end{array}
$$
Letting~$t\rightarrow\infty$, by the strong law of large numbers,~$(\ref{eq3.38})$ and~$(\ref{eq3.21})$ we deduce that
$$
\begin{array}{lll}
\displaystyle \limsup\limits_{t\rightarrow\infty}\frac{1}{t}\int^{t}_{0}\frac{1}{x(s)}\mathrm{d}s&\leq& \displaystyle \frac{1}{\sigma}\left(\limsup\limits_{t\rightarrow\infty}\frac{1}{t}\int^{t}_{0}\mu\psi(s)\mathrm{d}s+\delta+\frac{\sigma_{1}^{2}}{2}\right)\\
&=&  \displaystyle \frac{1}{\sigma}\Big[\frac{\mu}{\beta}\Big(\alpha-\frac{\sigma_{2}^{2}}{2}\Big)
+\delta+\frac{\sigma_{1}^{2}}{2}\Big].
\end{array}
$$
The proof is complete.
\qed
\begin{lemma}\label{le3.8}
 If~$\delta>h^2$ and~$\displaystyle
\alpha-\frac{\sigma_{2}^{2}}{2}-\frac{\sigma}{\delta-h^2}>0$, then the inequality
\begin{equation}\label{eq3.29}
\displaystyle \liminf\limits_{t\rightarrow\infty}\frac{1}{t}\int^{t}_{0}y(s)\mathrm{d}s\geq \lambda_{3}
\end{equation}
holds, where~$\lambda_{3}:= \displaystyle\frac{1}{\beta}\left(\alpha-\frac{\sigma_{2}^{2}}{2}-\frac{\sigma}{\delta-h^2}\right)$.
\end{lemma}
{\bf Proof.}
For any~$(x_0,y_0)\in\mathbb{R}_{+}^{2}$,  since $0<x(t)\leq \varphi(t)$ a.s.,  an application of the It\^{o} formula yields
$$
\begin{array}{lll}
\displaystyle \frac{1}{t}\ln y(t)&=&\displaystyle \alpha-\frac{\sigma_{2}^{2}}{2}-\frac{1}{t}\int^{t}_{0}\beta y(s)\mathrm{d}s-\frac{1}{t}\int^{t}_{0}x(s)\mathrm{d}s+\frac{\ln y_{0}}{t}+\frac{\sigma_{2}B_{2}(t)}{t}\\
&\geq&\displaystyle \alpha-\frac{\sigma_{2}^{2}}{2}-\frac{1}{t}\int^{t}_{0}\beta y(s)\mathrm{d}s-\frac{1}{t}\int^{t}_{0}\varphi(s)\mathrm{d}s+\frac{\ln y_{0}}{t}+\frac{\sigma_{2}B_{2}(t)}{t}.
\end{array}
$$
Taking ~$t\rightarrow\infty$, by~$(\ref{eq3.19})$ and~$(\ref{eq3.39})$ we have
$$
\begin{array}{lll}
\displaystyle \liminf\limits_{t\rightarrow\infty}\frac{1}{t}\int^{t}_{0}y(s)\mathrm{d}s&\geq& \displaystyle \frac{1}{\beta}\left(\alpha-\frac{\sigma_{2}^{2}}{2}-\liminf\limits_{t\rightarrow\infty}\frac{1}{t}\int^{t}_{0}\varphi(s)\mathrm{d}s\right)\\
&=&\displaystyle \frac{1}{\beta}\left(\alpha-\frac{\sigma_{2}^{2}}{2}-\frac{\sigma}{\delta-h^2}\right).
\end{array}
$$
The proof is complete.
\qed

Now, let us prove the existence of the invariant measure of the equation ~$(\ref{eq1.5})$.
\begin{theorem}\label{th3.6}
 If~$\delta>h^2$~and ~$\displaystyle\alpha-\frac{\sigma_{2}^{2}}{2}-\frac{\sigma}{\delta-h^2}>0$, then the process~$(x(t),y(t))$ has an invariant probability measure on~$\mathbb{R}_{+}^{2}$.
\end{theorem}
{\bf Proof.}
Let ~$\hbar$ and $H $ be two positive constants such that~$\hbar<\min\{\lambda_1, \lambda_2\},~H>\max\{\lambda_1, \lambda_2\}$, where $\lambda_1$ and $ \lambda_2$ defined in Theorem \ref{th3.5}  and Lemma \ref{le3.7}, respectively.
By  the H\"{o}lder inequality, we have
$$
\displaystyle \frac{1}{t}\int^{t}_{0}\mathbf{I}_{\{y(s)\geq\hbar\}}y(s)\mathrm{d}s\leq\displaystyle \left(\frac{1}{t}\int^{t}_{0}\mathbf{I}_{\{y(s)\geq\hbar\}}\mathrm{d}s\right)^{\frac{1}{2}}\left(\frac{1}{t}\int^{t}_{0}y^{2}(s)\mathrm{d}s\right)^{\frac{1}{2}},
$$
therefore
\begin{equation}\label{eq3.40}
\displaystyle \liminf\limits_{t\rightarrow\infty}\frac{1}{t}\int^{t}_{0}\mathbf{I}_{\{y(s)\geq\hbar\}}y(s)\mathrm{d}s\leq\displaystyle \left(\liminf\limits_{t\rightarrow\infty}\frac{1}{t}\int^{t}_{0}\mathbf{I}_{\{y(s)\geq\hbar\}}
\mathrm{d}s\right)^{\frac{1}{2}}\left(\limsup\limits_{t\rightarrow\infty}\frac{1}{t}\int^{t}_{0}y^{2}(s)\mathrm{d}s\right)^{\frac{1}{2}}.
\end{equation}
Moreover, we have
\begin{equation}\label{eq3.37}
y(t)\mathbf{I}_{\{y(t)\geq\hbar\}}=y(t)-y(t)\mathbf{I}_{\{y(t)<\hbar\}}\geq y(t)-\hbar.
\end{equation}
Hence, combing \eqref{eq3.40} with \eqref{eq3.37} yields
\begin{align}\label{eq3.30}
\begin{split}
\displaystyle\liminf\limits_{t\rightarrow\infty}\frac{1}{t}
\int^{t}_{0}\mathbf{I}_{\{y(s)\geq\hbar\}}\mathrm{d}s
\geq&\displaystyle \left(\liminf\limits_{t\rightarrow\infty}\frac{1}{t}\int^{t}_{0}\mathbf{I}_{\{y(s)\geq\hbar\}}
y(s)\mathrm{d}s\right)^{2} \left(\limsup\limits_{t\rightarrow\infty}\frac{1}{t}\int^{t}_{0}y^{2}(s)\mathrm{d}s\right)^{-1}\\
\geq&\displaystyle\left(\liminf\limits_{t\rightarrow\infty}
\frac{1}{t}\int^{t}_{0}(y(s)-\hbar)\mathrm{d}s\right)^{2} \left(\limsup\limits_{t\rightarrow\infty}\frac{1}{t}\int^{t}_{0}y^{2}(s)\mathrm{d}s\right)^{-1}\\
=&\displaystyle\left(\liminf\limits_{t\rightarrow\infty}\frac{1}{t}\int^{t}_{0}y(s)\mathrm{d}s
-\hbar\right)^{2} \left(\limsup\limits_{t\rightarrow\infty}\frac{1}{t}\int^{t}_{0}y^{2}(s)\mathrm{d}s\right)^{-1}\\
\geq&\displaystyle \frac{(\lambda_{3}-\hbar)^{2}}{\bar{M_{2}}}~~~~~a.s.
\end{split}
\end{align}
It follows from Lemma~$\ref{le3.7}$ that
\begin{equation}\label{eq3.31}
\displaystyle\limsup\limits_{t\rightarrow\infty}\frac{1}{t}\int^{t}_{0}\mathbf{I}_{\{x(s)\leq\hbar\}}\mathrm{d}s \leq  \hbar\limsup\limits_{t\rightarrow\infty}\frac{1}{t}\int^{t}_{0}\frac{1}{x(s)}\mathrm{d}s \leq \hbar \lambda_{2}~~~~~a.s.
\end{equation}
Similarly,  by Lemma~$\ref{le3.5}$ and $\ref{le3.6}$, one observes that
\begin{equation}\label{eq3.32}
\displaystyle\limsup\limits_{t\rightarrow\infty}\frac{1}{t}\int^{t}_{0}\mathbf{I}_{\{x(s)\geq H\}}\mathrm{d}s\leq\displaystyle\frac{1}{H}\limsup\limits_{t\rightarrow\infty}\frac{1}{t}\int^{t}_{0}x(s)\mathrm{d}s\leq\displaystyle\frac{M_{1}}{H}~~~~~a.s.
\end{equation}
and
\begin{equation}\label{eq3.33}
\displaystyle\limsup\limits_{t\rightarrow\infty}
\frac{1}{t}\int^{t}_{0}\mathbf{I}_{\{y(s)\geq H\}}\mathrm{d}s
\leq\displaystyle\frac{1}{H}\limsup\limits_{t\rightarrow\infty}\frac{1}{t}\int^{t}_{0}y(s)\mathrm{d}s\leq\displaystyle\frac{\bar{M_{1}}}{H}~~~~~a.s.
\end{equation}
Let~$A=\{(x,y):\hbar\leq x\leq H,~\hbar\leq y\leq H\}$. Choose more precise $\hbar$ and $ H $ such that $$\hbar<\displaystyle\min\left\{\lambda_1, \lambda_2, \frac{\lambda_{3}}{2},
~\frac{\lambda_{3}^{2}}{16\bar{M}_{2}\lambda_{2}}\right\}, ~H>\displaystyle\displaystyle\max\left\{\lambda_1,\lambda_2, \frac{16\bar{M}_{2}(M_{1}+\bar{M}_{1})}{\lambda_{3}^{2}}\right\}.
$$
From~$(\ref{eq3.30})$-$(\ref{eq3.33})$,  we derive
$$
\displaystyle\liminf\limits_{t\rightarrow\infty}\frac{1}{t}\int^{t}_{0}\mathbf{I}_{\{(x(s),y(s))\in A\}}\mathrm{d}s
\geq\displaystyle\frac{(\lambda_{3}-\hbar)^{2}}{\bar{M}_{2}}-\frac{M_{1}+\bar{M}_{1}}{H}-\hbar \lambda_{2}>\displaystyle\frac{\lambda_{3}^{2}}{8\bar{M}_{2}}~~~~~a.s.
$$
Taking expectation on both sides, we have
$$
\displaystyle\mathbb{E}\left[\liminf\limits_{t\rightarrow\infty}\frac{1}{t}\int^{t}_{0}\mathbf{I}_{\{(x(s),y(s))\in A\}}\right]\mathrm{d}s
>
\frac{\lambda_{3}^{2}}{8\bar{M}_{2}}.
$$
Using the Fatou lemma and the Fubini theorem yields
\begin{equation}\label{eq3.34}
\displaystyle\liminf\limits_{t\rightarrow\infty}\frac{1}{t}\int^{t}_{0}\mathbb{P}(s,(x,y),A)\mathrm{d}s>\displaystyle \frac{\lambda_{3}^{2}}{8\bar{M}_{2}},~~~~~\forall(x,y)\in\mathbb{R}_{+}^{2},
\end{equation}
where~$\mathbb{P}(t,(x,y),\cdot)$ is the transition probability of~$(x(t),y(t))$. Obviously, the Markov process~$(x(t),y(t))$ on the state space~$\mathbb{R}_{+}^{2}=\{x>0,y>0\}$  has the Feller property. Thus,~$(\ref{eq3.34})$ and {\cite{MT1993}} imply that $(x(t),y(t))$ has an invariant probability measure~$\pi_2(\cdot).$
\qed

We  now present  the uniqueness of the invariant measure of ~$(x(t),y(t))$.
\begin{theorem}\label{T6} Under the conditions of Theorem \ref{th3.6} the solution  $(x(t), y(t))$
of  \eqref{eq1.5} has a unique invariant measure.
\end{theorem}
{\bf Proof.}\,\, For convenience, let $ 2 \zeta:= \displaystyle  (\delta-h^2)\left(\alpha-\frac{\sigma_{2}^{2}}{2}\right)-  \sigma. $ Obviously, the given conditions imply $\zeta>0.$ Furthermore, choose a constant $c>0$ small sufficiently such that
\be \la{equ+1}c ({\delta+\sigma_1^2}) \leq    \sigma \zeta.
\ee
Define   $U: \RR_+^2 \rightarrow \RR_+$ by
\be\la{equ30}
U(x, y)= x+\frac{c}{x}  +y^2+\displaystyle(\delta-h^2)\ln(1+\frac{1}{y}).
\ee
Computing $ \Le  U (x,y)$ yields
\begin{eqnarray}
 &&  \Le  U (x,y) \nonumber\\
  & & \quad = \displaystyle
 \left( \sigma+\frac{\rho{x y }}{\eta+y }-\mu{x y }-\delta{x }\right)
  -c\left(\frac{\sigma}{x^2}+\frac{\rho{ y }}{x(\eta+y )}-\frac{\mu y}{x}-\frac{\delta+\sigma^2_1}{x}\right)+\left[(2\alpha +\sigma_2^2)y^2
  \right.  \nonumber\\
&&\qquad   \displaystyle \left.
-2\beta y^3-2xy^2 \right]
+(\delta-h^2)\left[-\frac{\sigma_2^2}{2(y+1)^2}-\frac{\alpha-\sigma_2^2-x}{ y+1 }+\frac{\beta y}{y+1}\right].
    \label{equ31}
\end{eqnarray}
Noting that the definition of $h$ in  \eqref{eq-f} and  $\delta>h^2, ~ \displaystyle  \frac{ \mu y}{x}\le \frac{ \sigma}{2x^2}+\frac{ \mu^2}{2\sigma}$, one observes that
\begin{eqnarray}
    \Le  U (x,y)
 & \leq &  \displaystyle
 [\sigma  -(\delta-h^2) {x }]
  -\frac{c \sigma}{2x^2}  +\frac{c(\delta+\sigma^2_1)}{x}+(2\alpha +\sigma_2^2+\frac{c\mu^2}{2\sigma})y^2
    \nonumber\\
&&   \displaystyle
-2\beta y^3
+(\delta-h^2)\left[-\frac{\sigma_2^2}{2(y+1)^2}-\frac{\alpha-\sigma_2^2-x}{ y+1 }+\frac{\beta y}{y+1}\right]
\label{equ32}
\end{eqnarray}
This yields that there exist positive constants $h, H, \tilde{h}, \tilde{H}$ such that 
\be\la{equ41}
    \Le  U (x,y)\leq -\zeta, ~~ ~~(x,y)\notin D: =\{(x, y): h<x<H,    \tilde{h}<y< \tilde{H}\}.
    \ee
By {\cite[pp.106-122]{K2012}},
$(x(t), y(t))$ is positive recurrent with respect to $D $. Then
the desired
assertion follows. \qed

Moreover, by \cite{IK1974} and \cite{K2012}, we have the following ergodicity result.

\begin{theorem}\label{th3.7}
 Under the conditions of  Theorem~\ref{th3.6}, the model~(\ref{eq1.5}) has a unique invariant probability measure~$\pi_2$ with support~$\mathbb{R}_{+}^{2}$. Moreover,
\begin{itemize}
\item[$(\mathrm{i})$] For any~$\pi_2$-integrable~$f(x,y):\mathbb{R}_{+}^{2}\rightarrow\mathbb{R}$, we have
$$
\displaystyle\lim\limits_{t\rightarrow\infty}\frac{1}{t}\int^{t}_{0}f(x(s),y(s))\mathrm{d}s=\displaystyle \int_{\mathbb{R}_+^2} f(x,y)\pi_2(\mathrm{d}x,\mathrm{d}y)~~~~~a.s.~~\forall(x(0),y(0))\in\mathbb{R}_{+}^{2}.
$$
\item[$(\mathrm{ii})$]  Let $\|\cdot\|_{\mbox{var}}$ denote  the total variation norm, for $(x,y)\in\mathbb{R}_{+}^{2}$, we have
$$
\displaystyle\lim\limits_{t\rightarrow\infty}\|\mathbb{P}(t,(x,y),\cdot)-\pi_2(\cdot)\|=0 ~~~~~\forall(x,y)\in\mathbb{R}_{+}^{2}.
$$
\item[$(\mathrm{iii})$] For any~$\varepsilon>0$, there is ~$\zeta\in(0,1)$ such that
$$
\displaystyle\liminf\limits_{t\rightarrow\infty}\mathbb{P}\Big(t,x,y,[\zeta,\zeta^{-1}]\times[\zeta,\zeta^{-1}]\Big)
>1-\varepsilon.
$$
\end{itemize}
\end{theorem}

For a biological system the property $(\mathrm{iii})$ of Theorem \ref{th3.6} is also called stochastic strong permanence.

\section{Examples and numerical simulations}

In this section, we mainly illustrate the effects of noise intensity on effector cells and tumor cells. We select the data in  \cite{KMT1994} and \cite{Siu1986} , see Table ~4.1 below.
\begin{center}
{\bf{Table~4.1: The Significance  and value of the parameters 
}}\\
\begin{tabular}{llllllll}
\hline
\bf{parameter}&\bf{Real value/unit}&\bf{Biological significance}\\
\hline
~~~$a$&0.18~/day&the intrinsic growth rate of the TC\\
~~~$b$&$2.0\times10^{-9}$ /day~& Reciprocal of environmental capacity of TC\\
~~~$s$&$1.3\times10^{4}$ /piece$\times$day ~& the normal rate of inflow into the tumor site for EC\\
~~~$d$&0.0412~/day ~&   the coefficient of destruction and migration of EC\\
~~~$g$&$2.019\times10^{7}$ piece ~&the positive constant in response functional\\
~~~$q$&0.1245~/day & ~$q=fK, ~\displaystyle K=\frac{k_1}{k_2+k_3+k_{-1}}$ \\
~~~$r_1$&$2.422\times10^{-10}$ /day$\times$piece& $r_1=Kk_3$\\
~~~$r_2$ &$1.101\times10^{-7}$ day$\times$piece& $r_2=Kk_2$\\
\hline
\end{tabular}\\
\end{center}
where~$f$~is the positive constant of response function, $k_1$ and $k_{-1}$~describe the rates of binding of EC to TC and detachment of EC from TC without damaging cells, $k_2$~is the rate at which
EC-TC interactions irreversibly program TC for lysis, and $k_3$~is the rate at which EC-TC interactions inactivate EC. The non-dimensional treatment of the equation is done by  selecting the order of magnitude scales $E_0$ and $T_0$ for the $E$ and $T$ cell populations, respectively,
where ~$E_0=T_0 = 10^6$ cells \cite{KMT1994}. Using the non-dimensionalization method in \cite{KMT1994} yields coefficients in the model (\ref{eq1.5})  as follows
\begin{equation*}
\begin{array}{llllll}
\displaystyle\sigma=\frac{s}{r_2E_{0}T_{0}}=0.1181,
&\rho=\displaystyle\frac{q}{r_2T_{0}}=1.131,
&\displaystyle\mu=\frac{r_1}{r_2}=0.00311,
&\!\!\!\!\!\!\! \displaystyle\delta=\frac{d}{r_2T_{0}}=0.3743,\\
\\
\displaystyle\alpha=\frac{a}{r_2T_{0}}=1.636,&\eta=\displaystyle\frac{g}{T_{0}}=20.19, & \beta=\displaystyle\frac{a b}{r_2}=3.272\times10^{-3}, &x_0=5,~~y_0=50.
\end{array}
\end{equation*}
In addition, applying the Milstein method in~Higham \cite{W}, we obtain the discrete equation as follows:
\begin{align}\label{eq4.1}
\!\! \left\{\!\!\!
\begin{array}{lll}
 \displaystyle x_{k+1}= x_{k}+ \Big(\sigma+\frac{\rho{x_ky_k}}{\eta+y_k}-\mu x_ky_k -\delta{x_k}\Big)\Delta t+\sigma_{1}{x_k}\sqrt{\Delta t}\xi_k+\frac{\sigma^2_{1}{x^2_k}}{2} \big(\xi^2_k-1\big)\Delta t,\!\!\!\!\\
 \displaystyle y_{k+1}=  y_k+ \Big(\alpha{y_k}-\beta{y^{2}_k}-x_ky_k\Big)\Delta t+\sigma_{2}{y_k}\sqrt{\Delta t}\eta_k+\frac{\sigma^2_{2}{y^2_k}}{2} \big(\eta^2_k-1\big)\Delta t ,\!\!\!\!
\end{array}
\right.
\end{align}
where~$\xi_k$, $\eta_k$ ($k=1,2,\ldots$) are two independent Gaussian random variables, and both obey the normal  distribution with mean ~0 and variance ~1.

\begin{expl}\label{ex4.1}
{\rm Choose the noise intensities~$\sigma_1=0.2$, $\sigma_2=2$ in the stochastic tumor-immune model~$(\ref{eq1.5})$.
Then we have
$$
2\alpha- \sigma_2^2=-0.728<0.
$$
Theorem $\ref{th3.5}$ tell us that the density of tumor cells ~$y(t)$ is exponentially decreasing, see the right side of Figure $\ref{fi4.1}$. Meanwhile, Theorem $\ref{th3.5}$ also shows that the distribution of effector cells $x(t)$ weakly converges to the unique invariant probability measure $\pi_1(\cdot)$, the inverse gamma distribution with  $a_3=19.715$ and $b_3=5.905$. To further illustrate the result of Theorem $\ref{th3.5}$, we use the $\mathrm{K}$-$\mathrm{S}$ test with a significance level of $0.05$ to check if the stationary distribution of $x^{-1 }(t)$ is the gamma distribution. At this level of significance, by $\mathrm{Matlab}$ we do confirm that the stationary distribution of $x^{-1}(t)$ is the Gamma distribution.
And because $x^{-1}(t)\sim G(19.715,~5.905)$ is equivalent to $x(t)\sim IG(19.715,~5.905)$, we know that the stationary distribution of $x(t) $ is the inverse gamma distribution.
\begin{figure}[htp]
  \begin{center}
\includegraphics[width=15cm,height=4.5cm]{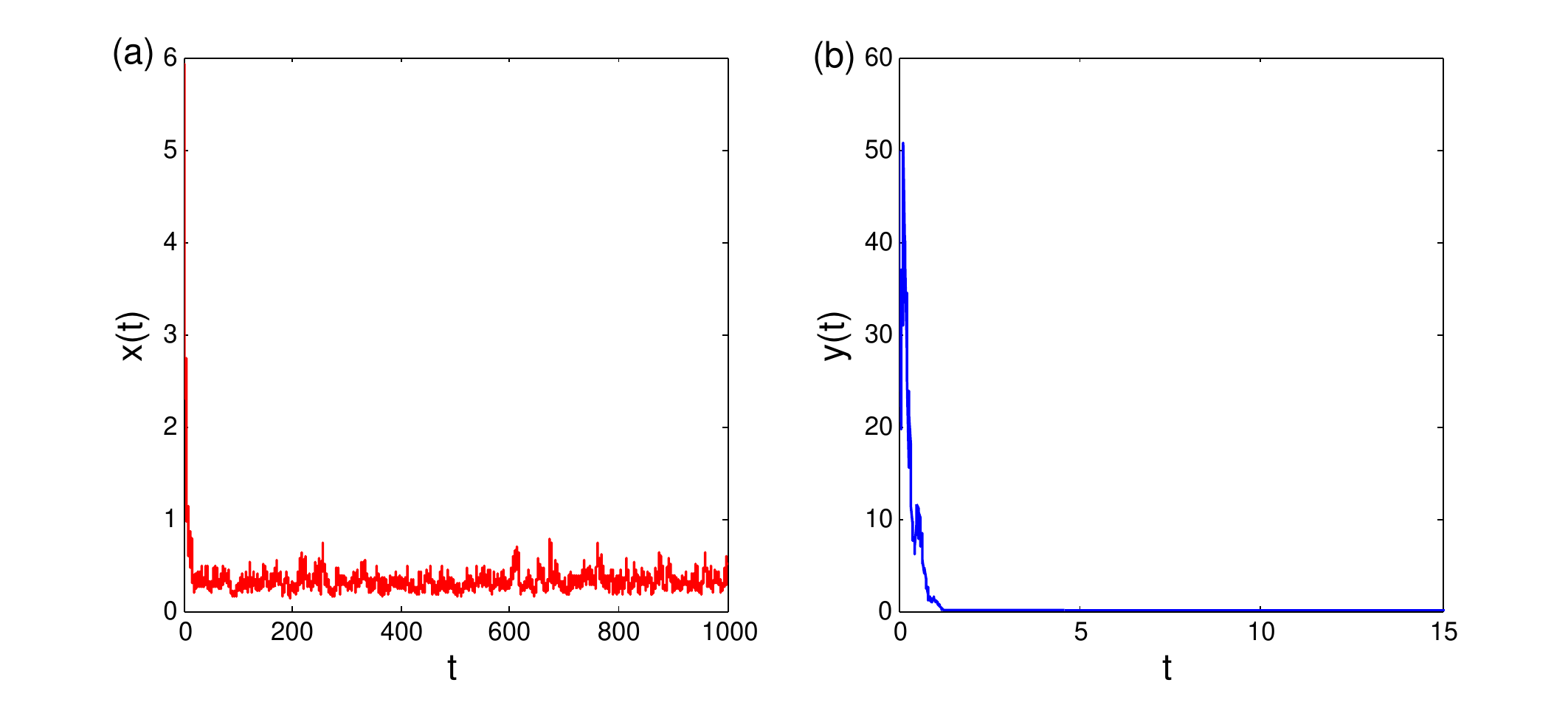}
   \end{center}
  \caption{Example~$\ref{ex4.1}$. For the stochastic tumor-immune model \eqref{eq1.5}, the red solid line depicts the density of effector cells ~$x(t)$, the solid blue line depicts the density of tumor cells ~$y(t)$.}
\label{fi4.1}
\end{figure}

\begin{figure}[htp]
  \begin{center}
\includegraphics[width=15cm,height=4.5cm]{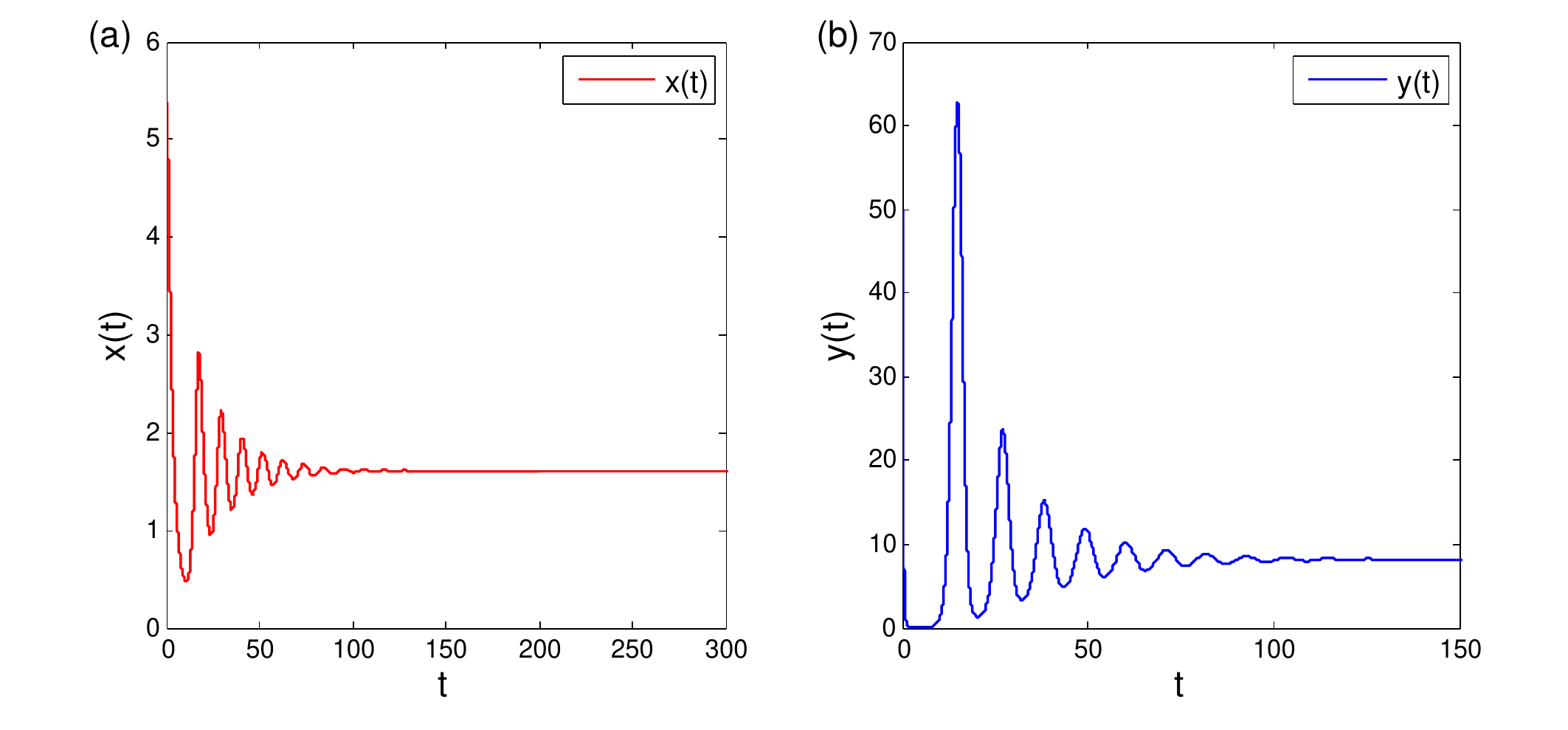}
   \end{center}
  \caption{Example~$\ref{ex4.1}$. For the deterministic tumor-immune model~$(\ref{eq1.4})$, the red solid line depicts the density of effector cells ~$x(t)$, the solid blue line depicts the density of tumor cells ~$y(t)$.}
  \label{fi4.1_2} \end{figure}

\begin{figure}[htp]
  \begin{center}
\includegraphics[width=10cm,height=5cm]{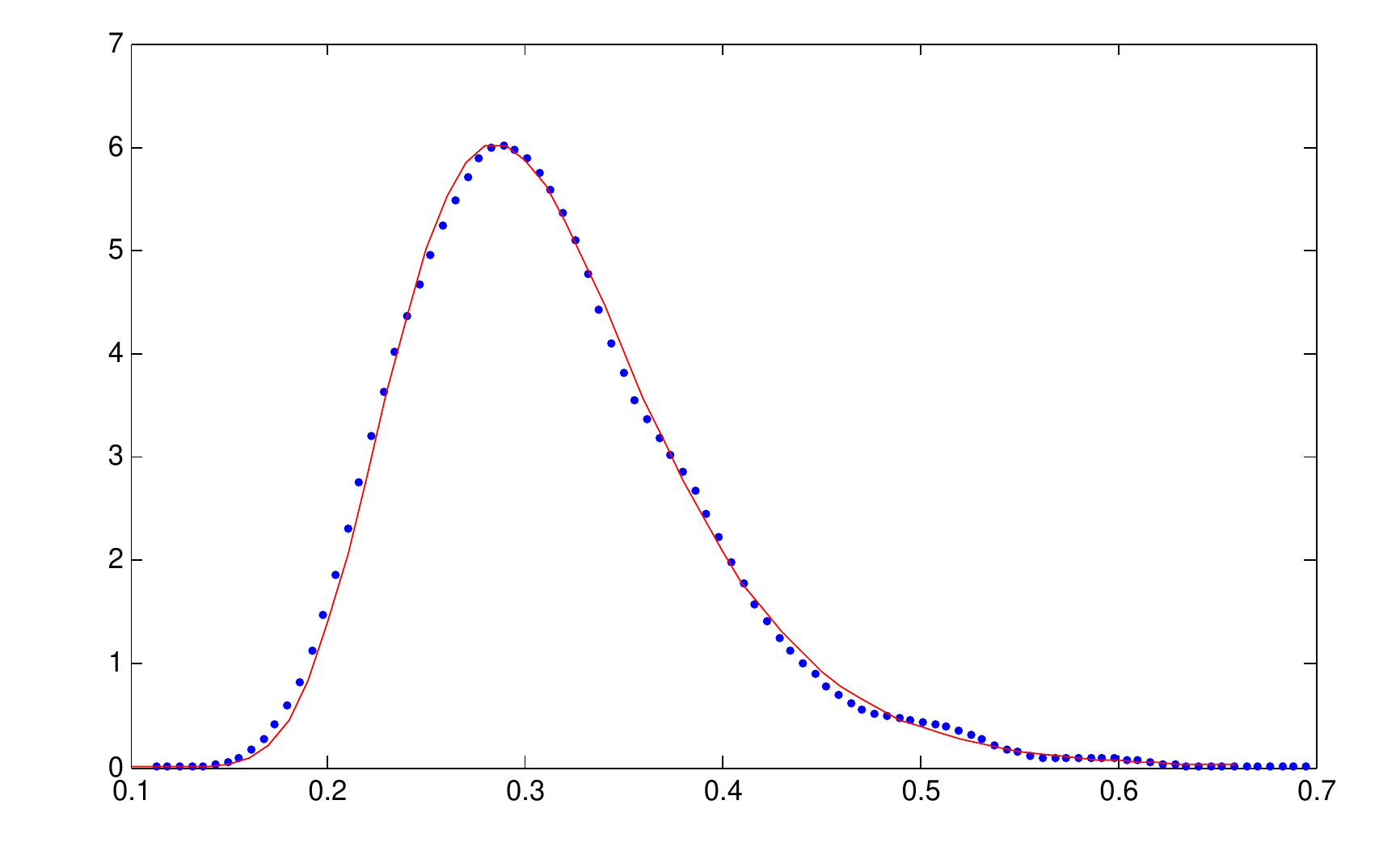}
   \end{center}
  \caption{Example~$\ref{ex4.1}$.
  The red solid line indicates the density function of  the inverse gamma~$IG(19.715,~5.905)$, the solid blue line indicates the empirical density function of the effector cell ~$x(t)$.}\label{fi4.2}
   \end{figure}
Furthermore, to more intuitively illustrate the result of Theorem $\ref{th3.5}$, we plot the empirical density function of $x(t)$ and the density function of the inverse gamma distribution $IG(19.715, ~5.905)$ in Figure $\ref{fi4.2}$. One observes obviously from the Figure $\ref{fi4.2}$ that as $t\rightarrow\infty$, the distribution of $x(t)$ weakly converges to $\pi_1(\cdot)$. Thus, this example illustrates the significance of the result of Theorem $\ref{th3.5}$. On the other hand, we  compare  the simulations of the stochastic model with the deterministic for the same parameter values. Figure $\ref {fi4.1_2}$ depicts that the path  of $y(t)$  in the  deterministic model  tends to a positive equilibrium with frequency vibration, namely, the tumor cells are not extinct. However, in Figure $\ref {fi4.1_2}$, one observes that when the noise intensity is large such that~$2\alpha<\sigma_2^2
 $, the tumor cells are extinct.  It is revealed that stochastic factors cannot be ignored, and their existence plays a key role in the permanence and extinction of the tumor cells.}
\end{expl}

\begin{expl}\label{ex4.2}
 {\rm In the stochastic tumor-immune model $(\ref{eq1.5})$, let $\sigma_2=0.25$, which implies that the stochastic environment has a weak effect on the intrinsic growth rate of tumor cells. At the same time, the binding rate of $\mathrm{EC}$ to $\mathrm{TC}$ will be decreased when the immune response of the effector cells to the tumor cells is weak or the tumor cells are less irritating to the effector cells. Therefore, in this example we reduce the binding rate $k_1$ of $\mathrm{EC}$ and $\mathrm{TC}$ in the literature \cite{KMT1994}, let $\rho=0.613$.
Compute
$$
\delta-h^2=0.09089>0,~~~~~ \alpha-\frac{\sigma^2_2}{2}-\frac{\sigma}{\delta-h^2}=0.30539>0.
$$
These imply that  the conditions of Theorems $\ref{th3.6}$ and $\ref{th3.7}$ hold. By virtue of Theorem $\ref{th3.6}$ and $\ref{th3.7}$ the solution $(x(t), y(t))$ of the stochastic tumor-immune model $(\ref{eq1.5})$ has a unique invariant probability measure $\pi^{* }$, and the system is stochastically  permanent.  Figure ~$\ref{fi4.3}$ depicts the trajectories of the effector cells $x(t)$ and the tumor cells $y(t)$ in $(\ref{eq1.5})$. Figure $\ref{fi4.4}$ is the phase diagram with respect to the model $(\ref{eq1.5})$. Figure $\ref{fi4.5}$ depicts the empirical density of the invariant measure $\pi^{*}$ of the stochastic model $(\ref{eq1.5})$.
Therefore, this example verifies the theoretical results of  Theorems ~$\ref{th3.6}$ and ~$\ref{th3.7}$ well. }
\begin{figure}[!htp]
  \centering
\includegraphics[width=15cm,height=5.5cm]{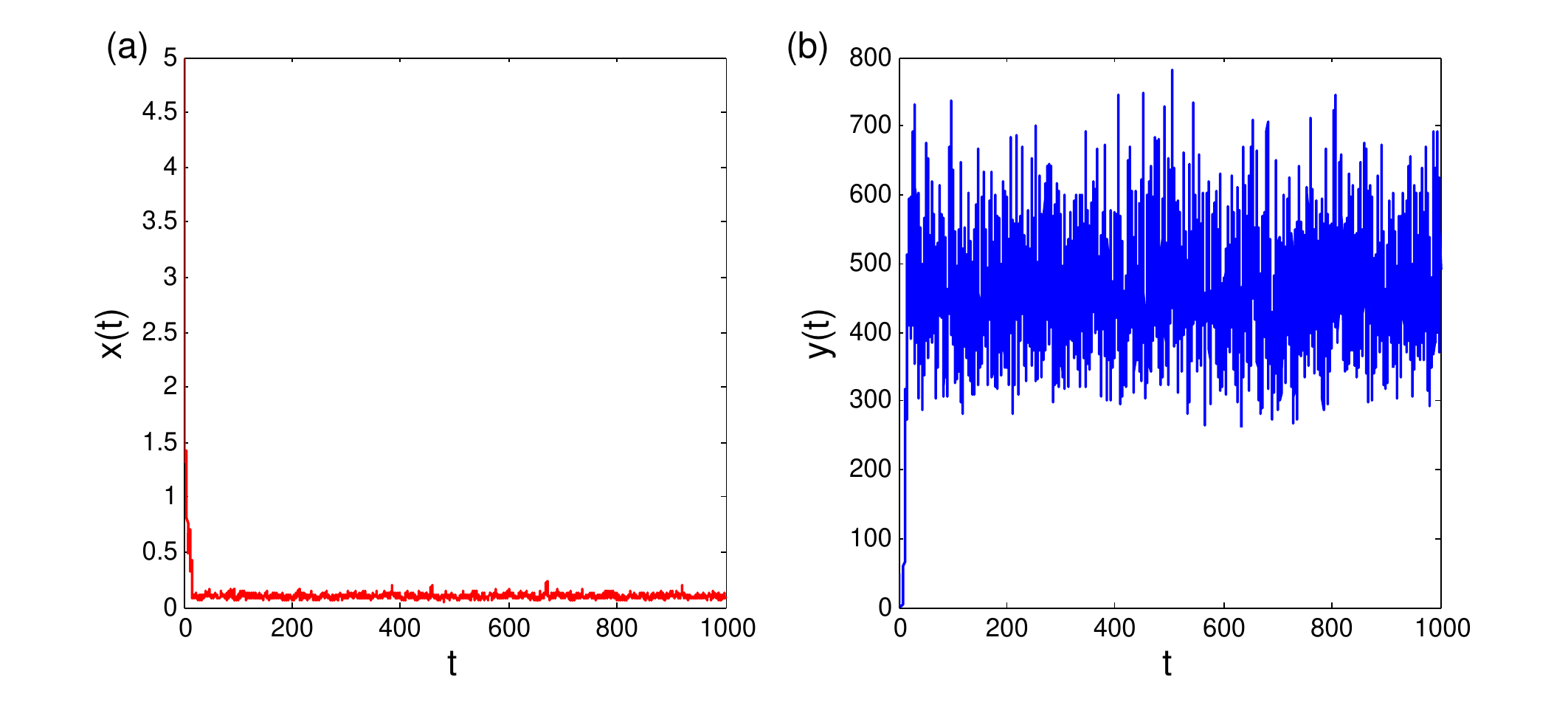}
  \caption{Example~$\ref{ex4.2}$. For stochastic tumor-immune model ~$(\ref{eq1.5})$, the red solid line indicates the effector cell ~$x(t)$, the solid blue line indicates the tumor cell ~$y(t)$.  }
  \label{fi4.3}
\end{figure}
\begin{figure}[htp]
  \centering
\includegraphics[width=12cm,height=6cm]{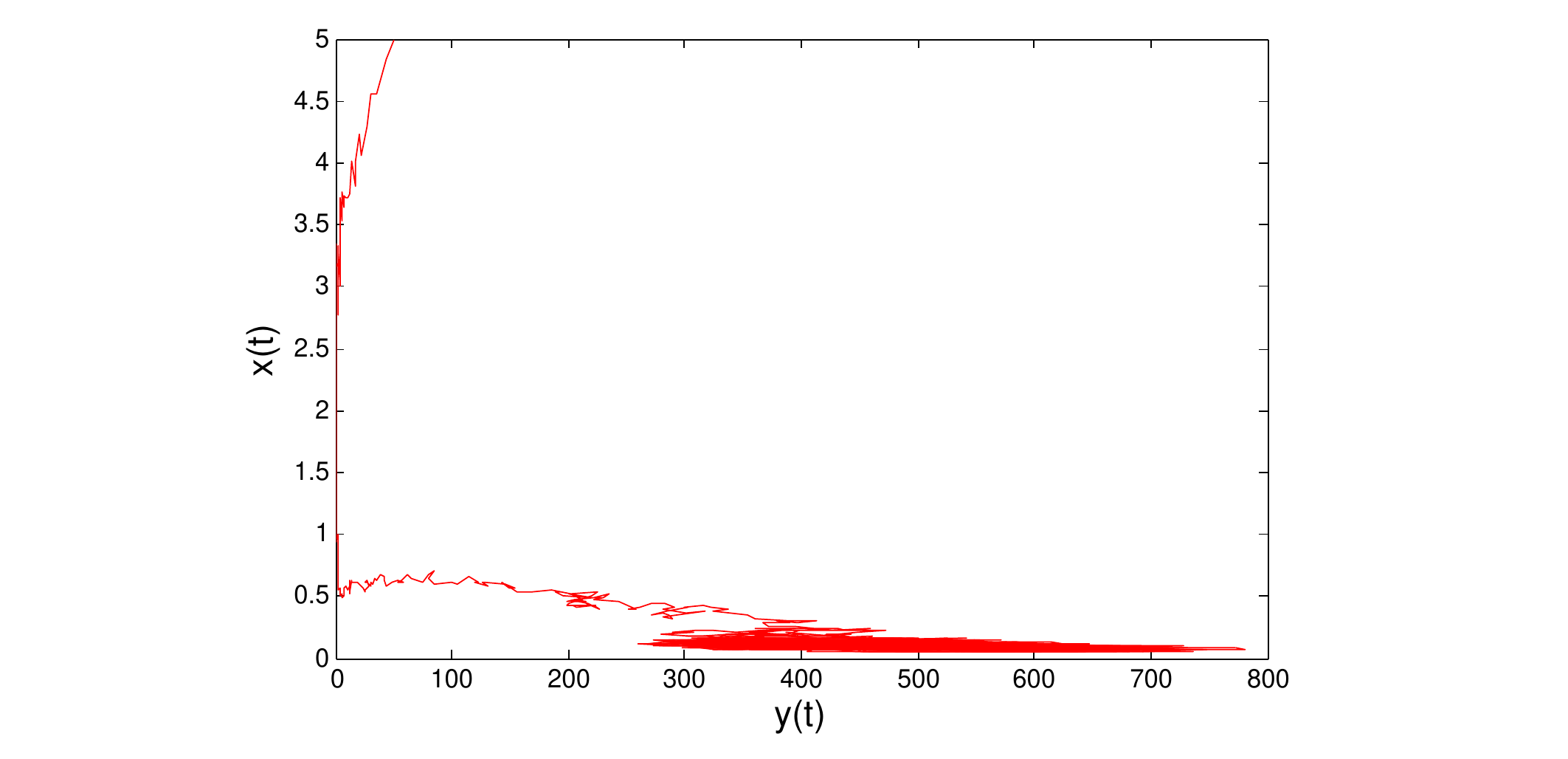}
  \caption{Example~$\ref{ex4.2}$. The phase diagram of the stochastic tumor-immune model~$(\ref{eq1.5})$.}
  \label{fi4.4}
\end{figure}
\begin{figure}[htp]
  \centering
\includegraphics[width=12cm,height=6.7cm]{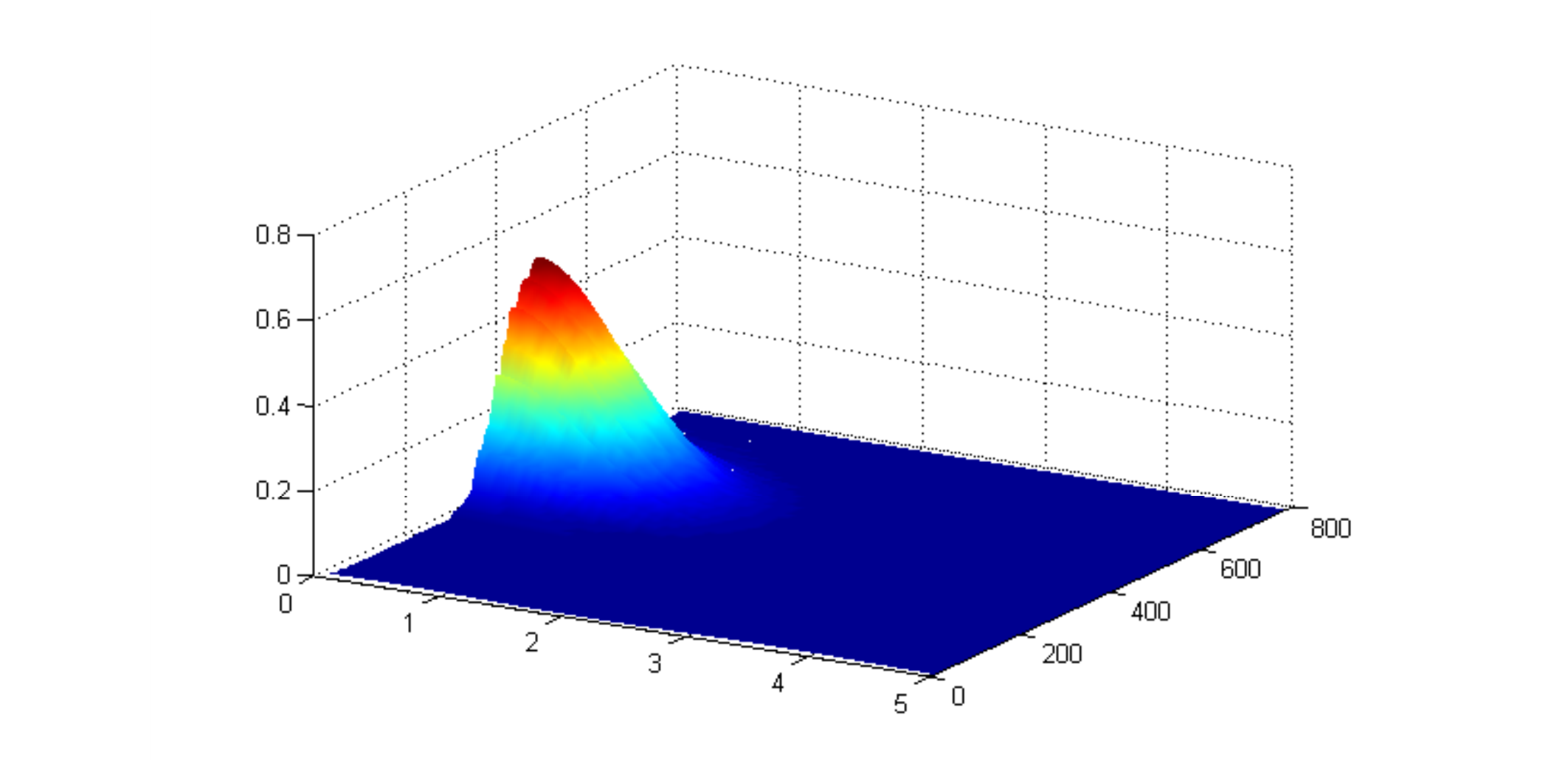}
  \caption{Example~$\ref{ex4.2}$.
   The empirical density of the invariant measure~$\pi^{*}$ of the stochastic tumor-immune model~$(\ref{eq1.5})$.
   }
  \label{fi4.5}
\end{figure}
\end{expl}

\section{Concluding remarks}

This paper mainly studies the dynamical  behaviors of the tumor-immune model proposed by  Kuznetsov and Taylor \cite{KMT1994} perturbed by the environment noise. Firstly, we prove the existence and uniqueness of the global positive solution for the tumor-immune system by the method  of stochastic Lyapunov analysis. Next, by constructing  appropriate comparison equations, we obtain the asymptotic moment boundedness of the effector cells and the tumor cells. Regarded the boundary equation $(\ref{eq3.24})$ as a bridge, it is pointed out that when tumor cells are subject to the strong noise, the density of the tumor cells decays to zero at an exponential rate while the density of effector cells tends to a stationary distribution.
 Furthermore, when the noise intensity of tumor cells is small relatively, by analyzing the upper and lower limits of the density of tumor cells and effector cells at time-average, we prove the existence and uniqueness of the stationary distribution of the stochastic tumor-immune model $(\ref{eq1.5})$. Moreover, the ergodicity and the stochastic permanence is  obtained . Finally, all of our main results are illustrated and  verified by numerical simulations. Overall, the fact is revealed  that the intensity of stochastic noise for tumor cells plays a key role in the permanence and extinction of tumor cells. As for strong noise,  compared with the deterministic, the dynamical behaviors of the stochastic model are different, and even richer.


\begin{thebibliography}{99}
{\small
\setlength{\baselineskip}{0.12in}
\parskip=0pt

 \bibitem{AB1997} Adam J A, Bellomo N. A Survey of Models for Tumor-Immune System Dynamics. Birkh$\ddot{\mathrm{a}}$user Boston, 1997.

  \bibitem{AG2006} Albano G, Giorno V. A stochastic model in tumor growth. Journal of Theoretical Biology, 2006, 242(2): 329-366.

 \bibitem{AM2004} Araujo R P, Mcelwain D L. A history of the study of solid tumour growth: the contribution of mathematical modeling. Bulletin of Mathematical Biology, 2004, 66(5): 1039-1091.

   \bibitem{H2}Bao J, Shao J. Permanence and extinction of regime-switching predator-prey models.  SIAM Journal on Applied Mathematics, 2016,  48(1): 725-739.

 \bibitem{Bellomo3} Bellomo N. Modeling complex living systems. A kinetic theory and stochastic game approach. Modeling and Simulation in Science, Engineering and Technology. Birkh$\ddot{\mathrm{a}}$user Bosten, 2008.

     \bibitem{BB2000} Burger R, Barton N H. The Mathematical Theory of Selection, Recombination, and Mutation. Wiley, New York, 2000.

\bibitem{DRW2005} DE Pillis L G, Radunskaya A E, Wiseman C L. A Validated Mathematical Model of Cell-Mediated Immune Response to Tumor Growth. Cancer Research, 2005, 65(17): 7950-7958.

  \bibitem{DND2015} Dieu N T, Nguyen D H, Du N H, et al. Classification of Asymptotic Behavior in A Stochastic SIR Model. Mathematics, 2015, 15(2): 1062-1084.

\bibitem{D2005} D'Onofrio A. A general framework for modeling tumor-immune system competition and immunotherapy: Mathematical analysis and biomedical inferences. Physica D Nonlinear Phenomena, 2005, 208(3-4): 220-235.

      \bibitem{DHY2015}Du N H, Hai N D, Yin G G. Conditions for permanence and ergodicity of certain stochastic predator-prey models. Journal of Applied Probability, 2015, 53(1): 187-202.

    \bibitem{EB2005} Elliott R L, Blobe G C. Role of Transforming Growth Factor Beta in Human Cancer. Journal of Clinical Oncology Official Journal of the American Society of Clinical Oncology, 2005, 23(9): 2078-2093.

      \bibitem{FBP2000} Ferrante L, Bompadre S, Possati L, et al. Parameter estimation in a Gompertzian stochastic model for tumor growth. Biometrics, 2000, 56(4): 1076-1081.

     \bibitem{GL1978} Garay R P, Lefever R. A kinetic approach to the immunology of cancer: stationary states properties of effector-target cell reactions. Journal of Theoretical Biology, 1978, 73(3): 417-438.

\bibitem{GR1974} Goel N S, Richter-Dyn N. Stochastic Models in Biology. Academic Press, 1974.

    \bibitem{W}  Higham D J. An algorithmic introduction to numerical simulation of stochastic differential equations. Society for Industrial and Applied Mathematics, 2001, 43(3): 525-546.

  \bibitem{IK1974} Ichihara K, Kunita H. A classification of the second order degenerate elliptic operators and its probabilistic characterization. Zeitschrift  f\"urWahrscheinlichkeitstheorie und Verwandte Gebiete, 1974, 30: 235-254.

 \bibitem{K2012} Khasminskii R. Stochastic Stability of Differential Equations. Springer Berlin Heidelberg, 2012.

\bibitem{Natalia}   Komarova N L , Wodarz D, Angelis E, et al. Selected Topics in Cancer Modeling: Genesis, Evolution, Immune Competition, and Therapy.  Birkh$\ddot{\mathrm{a}}$user Basel, 2008.

 \bibitem {KMT1994} Kuznetsov V A, Makalkin I A, Taylor M A, et al. Nonlinear dynamics of immunogenic tumors: Parameter estimation and global bifurcation analysis. Bulletin of Mathematical Biology, 1994, 56(2): 295-321.


\bibitem{LC2017} Li D, Cheng F. Threshold for extinction and survival in stochastic tumor immune system. Communications in Nonlinear Science and Numerical Simulations, 2017, 51: 1-12.

 \bibitem{MOE2016} Mahasa K J, Ouifki R, Eladdadi A, et al. Mathematical Model of Tumor-Immune Surveillance. Journal of Theoretical Biology, 2016, 404: 312-330.

   \bibitem {MAS2004} Mantovani A, Allavena P, Sica A. Tumour-associated macrophages as a prototypic type II polarised phagocyte population: role in tumour progression. European Journal of Cancer, 2004, 40(11): 1660-1667.


       \bibitem{AA} Mao X, Yuan C. Stochastic differential equations with Markovian switching. Imperial College Press, London, 2006.

      \bibitem{MT1993} Meyn S P, Tweedie R L. Stability of Markovian Processes III: Foster-Lyapunov criteria for continuous-time processes. Advances in Applied Probability, 1993, 25(3): 518-548.

     \bibitem{MB2009} Mukhopadhyay B, Bhattacharyya R. A Nonlinear Mathematical Model of Virus-Tumor-Immune System Interaction: Deterministic and Stochastic Analysis. Stochastic Analysis and Applications, 2009, 27(2): 409-429.


\bibitem{NM2000} Nowak M A, May R M. Virus dynamics: Mathematical Principles of Immunology and Virology. Oxford University Press, Oxford, UK, 2000.






\bibitem {OS1998} Owen M, Sherratt J. Modeling the macrophage invasion of tumors: Effects on growth and composition. Mathematics Applied in Medicine and Biology, 1998, 15: 165-185.

\bibitem{PG1997} Perelson A S, Ge R W. Immunology for physicists. Rev Mod Phys. Review of Modern Physics, 1997, 69(4): 1219-1268.



\bibitem{RDO2013} Riccardo C, Dumitru O, Oana C. Review of Stochastic Stability and Analysis Tumor-Immune Systems. Current Bioinformatics, 2013, 8(4): 390-440







\bibitem{Siu1986} Siu H, Vitetta E S, May R D, et al. Tumor dormancy. I. Regression of BCL1 tumor and induction of a dormant tumor state in mice chimeric at the major histocompatibility complex. Journal of Immunology, 1986, 137(4): 1376.




\bibitem{TC2006} Tan W Y, Chen C W. Cancer Stochastic Models. Encyclopedia of Statistical Sciences, 2006.


\bibitem{Villasana} Villasana M, Radunskaya A. A delay differential equation model for tumor growth. Journal of Mathematical Biology, 2003, 47(3): 270-294.

\bibitem{Wu}Wu J T, Kirn D H, Wein L M. Analysis of a three-way race between tumor growth, a replication-competent virus and an immune response. Bulletin of Mathematical Biology, 2004, 66(4): 605-625.

\bibitem{Yafia} Yafia R. Hopf bifurcation in differential equations with delay for tumor-immune system competition model. SIAM Journal on Applied Mathematics, 2007, 67(6):1693-1703.








    }

\end{thebibliography}
\end{document}